\documentclass[11pt]{amsart}
\usepackage{amsmath,amsfonts,latexsym,amssymb,amscd, graphicx,pictexwd}

\newcommand{\bpf}[1][Proof]{{\noindent {\sc #1: }}}
\newcommand{\epf}{{{\hfill $\Box$ \smallskip}}}
\newcommand{\ONE}{{\mathbf{1}}}
\newcommand{\NN}{{\mathbb N}}
\newcommand{\Fc}{{\mathcal F}}

\renewcommand{\P}{{\mathbb P}}
\newcommand{\E}{{\mathbb E}}
\newcommand{\R}{{\mathbb R}}

\newcommand{\ZZ}{{\mathbb Z}}

\newcommand{\HH}{{\mathbb H}}

\newcommand{\cal}{\mathcal}
\newcommand{\Co}{{\rm Co}}

\newcommand{\Z}{{\mathbb Z}}
\newcommand{\calF}{{\mathcal F}}
\newcommand{\calA}{{\mathcal A}}
\newcommand{\calC}{{\mathcal C}}

\newcommand{\eps}{{\varepsilon }}

\newcommand{\ceilt}{\lceil t\rceil}

\newtheorem{thm}{Theorem}[section]
\newtheorem{coro}[thm]{Corollary}
\newtheorem{lem}[thm]{Lemma}

\newtheorem{remark}[thm]{Remark}

\numberwithin{equation}{section}

\begin{document}
\title[Burgers Equation]{Space-time stationary solutions for the Burgers equation}

\author{Yuri Bakhtin}
\address{School of Mathematics\\ Georgia Institute of Technology\\ 686 Cherry Street\\ Atlanta, GA, 30332-0160, USA}
\email{bakhtin@math.gatech.edu}
\thanks{Yuri Bakhtin was supported by NSF CAREER Award DMS-0742424 and grant 040.11.264 from the Netherlands
Organisation for Scientific Research (NWO).  He is grateful for the hospitality of the Fields Institute in Toronto,
Delft Technical University, and CRM in Barcelona where parts of this work have been written.}

\author{Eric Cator}
\address{Radboud University Nijmegen\\
Heyendaalseweg 135\\ 6525 AJ  Nijmegen\\ The Netherlands}
\email{e.cator@math.ru.nl}
\thanks{Eric Cator is grateful for the hospitality of the Fields Institute in Toronto.}

\author{Konstantin Khanin}
\address{Department of Mathematics\\ University of Toronto \\ 40 St George Street\\
Toronto, Ontario, M5S 2E4, Canada}
\thanks{Konstantin Khanin is supported by NSERC Discovery Grant RGPIN 328565}
\email{khanin@math.toronto.edu}

\begin{abstract}

We construct space-time stationary solutions of the $1$D Burgers equation with random forcing in the absence of
periodicity or any other compactness assumptions. More precisely, for the forcing given by a homogeneous Poissonian
point
field in space-time we prove that there is a unique global solution with any prescribed average velocity. These global
solutions serve
as one-point random attractors for the infinite-dimensional dynamical system associated
to solutions to the Cauchy problem. The probability distribution of the global
solutions defines a stationary distribution for the corresponding Markov process.  We
describe a broad class of initial  Cauchy data for which the distribution of the
Markov process converges to the above stationary distribution.

\

Our construction of the global solutions is based on a study of the field of action-minimizing curves. We prove that for
an arbitrary value of the average velocity, with
probability 1 there exists a unique field of action-minimizing curves initiated at all
of the Poissonian points. Moreover action-minimizing curves corresponding to different starting points merge with each
other in finite time.

\end{abstract}


\maketitle

\section{Introduction}\label{sec:intro}

The Burgers equation is one of the basic hydrodynamic models. It describes the evolution of velocity field
of sticky particles that interact with each other only when they collide. In one space dimension, the inviscid
Burgers equation is
\begin{equation}
\label{eq:Burgers}
\partial_t u(x,t)+\partial_x\left(\frac{u^2(x,t)}{2}\right)=-\partial_x F(x,t).
\end{equation}
Here  $u(x,t)$ is the velocity of the particle located
at $x\in\R$ at time $t\in\R$. The external ``forcing'' term
$f(t,x)=-\partial_x F(x,t)$ describes the accelerations of particles. Although typically solutions of this equation
develop
discontinuities (shocks) in finite time, one can work with generalized solutions. So called entropy solutions, or
viscosity solutions, are globally well-defined and unique for a broad class of initial velocity profiles and forcing
terms.

\medskip

In this paper we study the long term behavior of the Burgers dynamics for the situation where
the forcing $f(x,t)$ is a space-time stationary random process.
In particular, we construct space-time stationary global solutions for the Burgers equation on the real line and show
that they can be viewed as one-point attractors.

The
results that we present here are connected with two big streams of research developed in the last twenty years.
The first one
concerns
stationary solutions and invariant measures for the randomly forced Burgers equation in a compact setting such as
periodic forcing that effectively reduces the system to dynamics on a circle, or torus in the
multidimensional version of equation~\eqref{eq:Burgers}, see
\cite{ekms:MR1779561}, \cite{Iturriaga:MR1952472}, \cite{Gomes-Iturriaga-Khanin:MR2241814}, or Burgers dynamics
 with random boundary conditions, see~\cite{yb:MR2299503}.

The main tool in the proof of these results is the Lax--Oleinik variational principle that
allows for an efficient analysis of the system via studying
the minimizers of the corresponding random Lagrangian system.
Namely, the velocity field can be represented as $u(t,x)=\partial_xU(t,x)$, where the potential $U(t,x)$
is a solution of the Hamilton--Jacobi equation
\begin{equation}
\label{eq:Hamilton-Jacobi-forced}
\partial_t U(x,t)+\frac{(\partial_x U(x,t))^2}{2}+F(x,t)=0.
\end{equation}
 The entropy solution to Cauchy problem for this equation
with initial data $U(\cdot,t_0)=U_0(\cdot)$ can be written as
\begin{equation}
\label{eq:LO}
U(x,t)=\inf_{\gamma:
[t_0,t]\to\R}\left\{U_0(\gamma(t_0))+{\int_{t_0}^t\left[\frac{1}{2}\dot\gamma^2(s)-F(\gamma(s),s)\right]} ds\right\},
\end{equation}
where the infimum is taken over all absolutely continuous curves $\gamma$
satisfying $\gamma(t)=x$.

The long-term behavior of the Burgers dynamics thus can be understood by studying this variational problem over long
time intervals. For example, to construct stationary solutions, one has to study the behavior of Lagrangian action
minimizers on an interval $[-T,t]$ as $T\to\infty$.  In particular, it was
proved in \cite{ekms:MR1779561} for the one-dimensional case and in \cite{Iturriaga:MR1952472} for the multi-dimensional
case
that in the compact situation, for each space point~$x$, the finite-time minimizers on $[-T,t]$ with endpoint $x$
stabilize to a unique infinite one-sided minimizer. Moreover, there is a unique
global minimizer which is a hyperbolic trajectory of the stochastic Lagrangian flow. The uniqueness of one-sided
minimizers implies the following One Force -- One Solution (1F1S) principle  for the Burgers equation: for almost
every realization of the forcing, there is a unique global solution with given average velocity, and at any given time
it
depends only on the history of the forcing up to that time.

The compactness condition in
the
above papers was
extremely important. Very few results have been available in noncompact setting. These results were dealing with the
case where the forcing was spatially nonhomogeneous and effectively compact, see~\cite{Khanin-Hoang:MR1975784}
and~\cite{Bakhtin-quasicompact}. The aim of our paper is to consider
forcing with homogeneous probability distribution in space and time. In other words, our goal is to replace the exact
periodicity (compactness) with
statistical homogeneity, namely, translational invariance in distribution.

Although the 1F1S principle is a simple conceptual fact in the compact situation, it is far from obvious that
stabilization of solutions occurs in the noncompact setting. Indeed, one can easily imagine a scenario where the
solutions
corresponding to different $T$ will be mostly influenced by the realization of the forcing in completely different
spatio-temporal patches. The main challenge in the problem,  is to show that this is not the case.
To tackle this difficulty, we use the second stream of research mentioned above, the theory of one-sided
infinite geodesics for first passage and last passage percolation models.

A simple model of that kind can be introduced in the following
way: consider a random potential on the lattice $\Z^2$. For a finite path on the lattice its action is
defined as the sum of all values of the potential along the path. A path is called a geodesic
between
two points if it minimizes the action among all paths connecting these two points.
An infinite path is called an infinite geodesic if every finite part of it is a geodesic between its endpoints.
The theory developed by Howard and Newman in \cite{Ne},\cite{HoNe3},\cite{HoNe2},\cite{HoNe}
describes the structure of such infinite length geodesics. It turns out that
each infinite geodesic has an asymptotic direction and for a fixed direction the geodesic
\emph{tree} spans the entire lattice. These type of results can be extended to a Poissonian
point field
setting where the disorder is due to random distances between configuration points. The problem is also closely related
to the study of
 Hammersley process (see
\cite{Aldous-Diaconis:MR1355056},\cite{Johansson:MR1757595},\cite{Wu},\cite{CaPi},\cite{CaPi-ptrf}) which is
based on optimal upright paths through Poissonian clouds. In all cases the existence of infinite geodesics is
proved using techniques going back to the work of Kesten \cite{Kesten:MR1221154}.

In the present paper we use  in a systematic way both Lagrangian methods developed in Aubry--Mather theory for
action-minimizing trajectories and weak KAM theory for Hamilton--Jacobi
equation (see \cite{Fathi},\cite{E:MR1682812}), and probabilistic techniques related to the first passage percolation
problem. We consider
equation~\eqref{eq:Burgers} and assume that the random forcing $f(x,t)=-\partial_xF(x,t)$ is
associated with a space-time homogeneous Poisson point field. One can think of every Poissonian point as a source of
a localized delta-type potential. This means that
$\int_{t_0}^t F(\gamma(s),s)ds$ in \eqref{eq:LO} is equal to the number of
Poissonian points visited by the curve $\gamma(s), s \in [t_0,t]$. Such singular random potentials were first
introduced
in the context of Burgers equation in \cite{Bakhtin-quasicompact}.

For systems with this type of forcing, shocks
are created at the Poissonian points. One can say that the main focus of
analysis is to study forward
dynamics of shocks and their merging and the backward dynamics of minimizers (action-minimizing curves) and their
coalescence. The crucial difference between \cite{Bakhtin-quasicompact} and the present paper is that the
Poissonian field in \cite{Bakhtin-quasicompact} is assumed to be
very nonhomogeneous in space with density decaying to zero at infinity. This results in quasi-compactness for the
problem, and the asymptotic analysis of the model is related to the one in the compact case. In particular there exists
a unique global minimizer and all other minimizers coalesce with it (see also~\cite{Khanin-Hoang:MR1975784}).
By imposing the homogeneity condition on the density of the Poissonian point field, we are forced to deal with a
completely noncompact situation. Despite the fact that the global behavior here is quite different (in particular, a
global minimizer does not exist), one still can prove almost sure existence and uniqueness of global solutions for fixed
``average'' velocity. The global solutions are stationary in time and space which reflects the translation-invariance
of the probability distribution for the forcing. The construction is based on a study of a global tree of
action-minimizing curves associated with every Poissonian point.



\medskip

We finish the introduction with a discussion of future directions. We believe that the
theory we discuss in this paper opens a wide
research area and below we formulate the most interesting and important open problems.
We plan to address some of the problems listed below in the future.

The first step in developing a general theory is elimination of the singular character of the forcing. There are
several settings where such elimination looks plausible. Consider a smooth nonnegative potential $F(x,t)$
vanishing outside of a disk of small radius $\eps>0$ and equal to $\eps^{-1}$ everywhere in the disk except for a thin
boundary layer. We then consider a homogeneous Poissonian point field and assume that every Poissonian point $(x_i,t_i)$
contributes $F(x-x_i,t-t_i)$ to the external potential, so that the total potential is given by the sum of these
contributions over all configuration points. In a certain sense the system considered in the present paper corresponds
to the limit $\eps\to 0$. We believe that similar results hold for the model with small positive  $\eps>0$. Namely, for
$\eps$ small enough, there exists a stationary global solution that approaches the global solution constructed in our
paper as $\eps\to 0$.

Another interesting generalization concerns systems with kick forcing. The corresponding variational problem will look
the following way.
\[
 U_v(x_0,0)=\min_{(x_i)_{i=-\infty}^0}\left[
\sum_{i}A_i(x_i,x_{i+1})+\frac{(x_{i+1}-x_i-v)^2}{2}\right].
\]
 Here the right-hand side can be defined up to an additive constant. This means that the increments
$U(x_0,0)-U(0,0)$ are well-defined provided that minimizing orbits for different starting points $x_0$ are
asymptotic to each other. We shall assume that random kernels $A_i$ are independent copies of a double stationary
process $A(x,y)$ (i.e., the distribution of the process $A(x+a,y+b)$ does not depend on $a,b$) with fast
enough decay of correlations.
The condition of double stationarity is important since it guarantees that the shape function (see
Section~\ref{sec:subadd}) is quadratic. However, it is natural to conjecture that existence and uniqueness of a global
solution only requires that the distribution of  $A_i(x+a,y+a)$ does not depend on $a$.

Of course, the most general result  will apply to any reasonable stationary forcing potential $F(x,t)$ in
equation \eqref{eq:Hamilton-Jacobi-forced} with fast decay of correlations. However, this problem is technically too
hard at present. It is also
not clear how far this program can be pushed in higher dimensions. We believe that
extension of our results in this direction will use Lagrangian methods in a very essential way.

Another set of problems is connected to the viscous case. In the case of singular potential corresponding to Poissonian
points, positive viscosity makes them irrelevant. However, for the smoothed problems discussed above, the
interplay between the viscosity and potential is nontrivial. It will be interesting to construct global stationary
solutions for positive viscosity $\nu$ and show that they approach the inviscid ones as $\nu\to 0$. In the case of
positive viscosity, the solution on a finite time
interval is determined not by a single action minimizing path, but by a random probability distribution on the space
of paths of fixed length
(directed polymers). It is tempting to make a conjecture that with probability 1 these random probability
distributions have a random limit as the time interval $[-T,t]$ converges to
$(-\infty, t]$ and the endpoint of the paths $\gamma(t)$ stays fixed. The limiting
probability distributions can be considered as Gibbs measures on the interval
$(-\infty,t]$. These Gibbs measures most probably cannot be extended to $(-\infty, +\infty)$.
The obstruction to such an extension is related to strong fluctuations of Gibbs
measure as $t\to \infty$. This is similar to the mechanism preventing existence
of global minimizers.

\medskip

Finally we briefly discuss connection with a very active area of KPZ scalings.
The global solution to the random Hamilton--Jacobi equation
$U(x,t)$ constructed in the
present paper is a process with stationary increments. If the increments are also
weakly dependent for distant intervals then one should expect to have a Gaussian behavior  with usual CLT scalings for
increments on large intervals. Namely, for some $\sigma>0$,
\begin{equation}
\label{eq:g}
\frac{U(x,0) - U(0,0) }{\sqrt{|x|}} \stackrel{distr}{\longrightarrow} \mathcal{N}(0,\sigma^2).
\end{equation}
In fact, an even stronger statement of the type of the invariance principle should
hold true. It is quite easy to see that such a ``non-exotic'' asymptotic behavior
will imply KPZ scalings for the limiting minimizer. To make a precise statement,
denote by $\gamma(t), t \in (-\infty,0]$ the one-sided minimizer with
zero average velocity which
originates at the origin at time zero. Then fluctuations of $\gamma(-t)$
and $U(\gamma(-t), -t) - U(0, 0)$ must be of the order $t^{2/3}$ and $t^{1/3}$ respectively. This
follows from (\ref {eq:g})
and the balance condition between $U(\gamma(-t), -t) - U(0, -t)$ and the contribution coming from
the shape function.
It is expected that the probability distribution for
$$\left(\frac{\gamma(-t)}{t^{2/3}}, \, \frac{U(\gamma(-t), -t) - U(0, 0)}{t^{1/3}}\right)$$
converges to the universal limit which is related to Tracy--Widom distribution
for the GOE random matrix ensemble. The exact form of the limiting distribution cannot be derived from the qualitative
argument above. However, the
universality of the limit law may be a more realistic target. We believe that
the limiting behavior is determined by the asymptotic properties of stationary
random point field of shocks, i.e., the points where
$U(x,0)$ is not smooth. Notice, however, that the field of shocks at a given time is
equipped with a random parameter, namely age of a shock, attached to
every point. The age is a time interval for which a given shock can be traced in
the past. It looks plausible that asymptotic statistical properties of the
point field of ``aged'' shock  determines the KPZ-type limits completely.

\section{The setting and notation}\label{sec:notation_and_results}

The forcing in our system is given by a Poissonian point field $\omega$ on space-time $\R\times\R=\R^2$ with Lebesgue
intensity measure. Throughout the paper we adopt the picture where the space axis of $\R^2$ is horizontal and
the time axis is vertical and directed upward.

The configuration space $\Omega_0$ is the space of all locally
finite point sets in space-time. For a Borel set $A\subset\R^2$,
we use $\omega(A)$ to denote the number of Poissonian points in $A$, and introduce the $\sigma$-algebra $\Fc_0$ generated by maps $\omega\mapsto
\omega(A)$ with  $A$ running through all bounded Borel sets. The probability measure $\P_0$ is such that
for any bounded Borel set $A$, $\omega(A)$ is a Poisson random variable with mean equal to the Lebesgue measure of $A$,
and for
disjoint bounded Borel sets $A_1,\ldots,A_m$, the random variables $\omega(A_1),\ldots,\omega(A_m)$ are independent.

Often, we treat the point configuration $\omega$ as a
locally bounded Borel measure with a unit atom at each point of the configuration.
For background on Poisson point fields, also called Poisson processes, we refer to
\cite{Daley:MR1950431}.

There is a natural family of time-shift operators $(\theta^t)_{t\in\R}$ on the Poisson configurations: the configuration
$\theta^t\omega$ is obtained from $\omega$ by shifting each point $(x,s)\in\omega$  to $(x,s-t)$.

The space of velocity potentials that we will consider will be~$\HH$, the space of all locally
Lipschitz functions $W:\R\to\R$ satisfying
\begin{align*}
 \liminf_{x\to+\infty}\frac{W(x)}{x}&>-\infty,\\
 \limsup_{x\to-\infty}\frac{W(x)}{x}&<+\infty.\\
\end{align*}
Although it is possible to work with weaker conditions, some restrictions on the growth rate
of $W(x)$ as $x\to\pm\infty$ are necessary to control velocities of particles coming
from $\pm\infty$.

Let us define random Hamilton--Jacobi--Hopf--Lax--Oleinik (HJHLO) dynamics on~$\HH$.
For a function $W\in\HH$, a Poissonian configuration $\omega$, and an absolutely continuous trajectory (path) $\gamma$
defined on $[s,t]$, we introduce the action
\begin{equation}
\label{eq:action}
A_\omega^{s,t}(W,\gamma)=W(\gamma(s))+S^{s,t}(\gamma)-\omega^{s,t}(\gamma),
\end{equation}
where
\[
 S^{s,t}(\gamma)=\frac{1}{2}\int_s^t\dot\gamma^2(r)dr
\]
is the kinetic action, and
$\omega^{s,t}(\gamma)=\omega(\{(\gamma(r),r):\ r\in[s,t)\})$ denotes the number of configuration points that $\gamma$
passes through. The last term in~\eqref{eq:action} is responsible for the interaction with the external forcing
potential corresponding to the realization of the Poissonian field.

We now consider the following minimization problem:
\begin{align}
 \label{eq:optimization_problem}
 A_\omega^{s,t}(W,\gamma)&\to \inf,\\
\gamma(t)&=x.\notag
\end{align}

Notice that the optimal trajectories are given by straight lines for any time interval on which the trajectory stays
away from the configuration points. Since  Poissonian configurations
are locally finite, it is sufficient to take the minimum over broken lines with
vertices at configuration points.


\begin{lem}\label{lem:full_measure_set_where_dynamics_is_defined} There is a set $\Omega_1\in\Fc_0$ with
$\P(\Omega_1)=1$
such that for all $t\in\R$
$\theta^t\Omega_1=\Omega_1$, and
for any $\omega\in\Omega_1$,
 for any $W\in\HH$, any $x\in\R$ and any $s,t$ with $s<t$ there is a path $\gamma^*$ that realizes
the minimum in~\eqref{eq:optimization_problem}. The path $\gamma^*$ is a broken line with finitely many segments, all
its
vertices belong to $\omega$.
\end{lem}
The proofs of
this and other statements of this section are
given in Section~\ref{sec:proofs_of_basic_burgers_facts}.

We denote the restrictions of $\Fc_0$ and $\P_0$ onto $\Omega_1$ by $\Fc_1$ and $\P_1$.
Since this restriction $\P_1$ still defines a Poisson point field with the same intensity measure, from
now on for convenience we remove from $\Omega_0$ the zero measure complement to $\Omega_1$ and work with the probability
space $(\Omega_1,\Fc_1,\P_1)$.

We denote the
infimum (minimum) value in~\eqref{eq:optimization_problem} by $\Phi^{s,t}_\omega W(x)$. The family of random nonlinear
operators
$(\Phi^{s,t}_\omega)_{s\le t}$ is the main object in this paper. Our main goal is to understand the
asymptotics of $\Phi^{s,t}_\omega$ as $t-s\to\infty$.

We will need several properties of the
random nonlinear operator $\Phi^{s,t}_\omega$ defined on $\HH$.
We begin with a lemma that shows that $\Phi^{s,t}_\omega W$ can be understood as the potential of the velocity field
given by the terminal velocities of minimizers and that it
produces a tessellation of space-time into the domains of influence of configuration points.

For a Borel set $B\subset \R^2$, we denote the restriction of
$\omega$ to $B$ by $\omega\bigr|_B$.

\begin{lem}\label{lem:evolution_on_Burgers_potentials} For any $\omega\in\Omega_1$, $W\in\HH$,
$s,t\in\R$ with $s<t$, the following holds true: \begin{enumerate}
 \item \label{it:open-domains-of-influence}  For any $p\in\omega|_{\R\times[s,t)}$, the set $O_p$ of points $x\in\R$
such that $p$ is the last
configuration point visited by a unique minimizer for problem~\eqref{eq:optimization_problem}, is open. Also, the set of
points with a unique minimizer that does not pass through any configuration
points is open. The union of these open sets is dense in $\R$.
\item \label{it:two-ways-to-compute-velocity}
If $x_0$ belongs to
one of these open sets, $\gamma(t)=x_0$, and $A^{s,t}_\omega(W,\gamma)=\Phi_\omega^{s,t}W(x_0)$, then
$\Phi_\omega^{s,t}W(x)$ is
differentiable at
$x_0$ w.r.t.\ $x$, and
\[
\frac{d}{dx}\Phi_\omega^{s,t}W(x)\bigr|_{x=x_0}=\dot\gamma(t).
\]
At a boundary point $x_0$ of any of the open sets introduced above, the right and left derivatives of
$\Phi_\omega^{s,t}W(x)$ w.r.t.~$x$
are well defined. They are equal to the slope of, respectively, the leftmost and rightmost minimizers realizing
$\Phi_\omega^{s,t}W(x_0)$.
\item \label{it:piecewise_linearity} For any $p\in\omega|_{\R\times[s,t)}$, the function $x\mapsto
\frac{d}{dx}\Phi_\omega^{s,t}W(x)$ is linear
in~$O_p$.
\item \label{it:continuity_of_HJ_solution}The function $x\mapsto\Phi_\omega^{s,t}W(x)$ is locally Lipschitz.

\end{enumerate}
\end{lem}

The following statement is the cocycle property for the operator family~$(\Phi^{s,t}_\omega)$. It is a direct
consequence of Bellman's principle of dynamic programming.

\begin{lem}\label{lem:cocycle} If $\omega\in\Omega_1$, then for any $W\in\HH$, any $s,r,t$ satisfying
$s<r<t$,
$\Phi_\omega^{r,t}\Phi_\omega^{s,r}W$ is well-defined and equals $\Phi_\omega^{s,t}W$. If $\gamma$ is an optimal path
realizing $\Phi_\omega^{s,t}W(x)$, then the restrictions of $\gamma$ on $[s,r]$ and $[r,t]$
are optimal paths realizing $\Phi_\omega^{s,r}W(\gamma(r))$ and   $\Phi_\omega^{r,t}(\Phi_\omega^{s,r}W)(x)$.
\end{lem}

Introducing $\Phi^t_{\omega}=\Phi^{0,t}_{\omega}$ we can rewrite the cocycle property as
\[
 \Phi^{0,t+s}_{\omega}W=\Phi^{s}_{\theta^t\omega} \Phi^{t}_\omega W,\quad s,t> 0,\quad\omega\in\Omega_1.
\]

Note that potentials are naturally defined up to an additive constant. It is thus convenient to
work with $\hat \HH$, the space of equivalence classes of potentials from $\HH$. The cocycle $\Phi$ can be
projected on $\hat\HH$ in a natural way. We denote the resulting cocycle on $\hat\HH$ by $\hat\Phi$.

Let us now explain how the dynamics that we consider is connected to the classical Burgers
equation. One way to describe this connection is to introduce a
mollification of the Poisson integer-valued measure. We give the corresponding statement without a proof. Let us take
smooth kernels $\phi,\psi:\R\to[0,\infty)$ with bounded support, satisfying $\int_{\R}\phi(t)dt=1$ and
$\max_{x\in\R}\psi(x)=1$, and for each $\eps>0$ consider the potential of shot-noise type:
\[
F_\eps(x,t)=-\frac{1}{\eps}\sum_{(y,s)\in\omega}\phi\left(\frac{t-s}{\eps}\right)\psi\left(\frac{x-y}{\alpha(\eps)}
\right),
\]
where $\alpha$ is any function satisfying $\lim_{\eps\downarrow 0}\alpha(\eps)=0$.

\begin{lem} With probability 1, for all $s,t,x\in\R$ and $W\in\HH$, the entropy solution $U_\eps(x,t)$ of the Cauchy
problem for the
Hamilton--Jacobi equation with smooth forcing potential $F_\eps(\cdot,\cdot)$ converges, as $\eps\to 0$, to
$U(x,t)=\Phi^{s,t}_\omega W(x)$.
\end{lem}

The next statement shows that
away from the Poissonian points the system we consider behaves as unforced Burgers dynamics.
\begin{lem}\label{lem:Hamilton-Jacobi} For
all $\omega\in\Omega_1$, $s\in\R$, $W\in\HH$, the function
$U(x,t)=\Phi^{s,t}_\omega W(x)$ is an entropy solution of the Hamilton--Jacobi equation
\begin{equation}
\label{eq:Hamilton-Jacobi}
\partial_t U+\frac{(\partial_x U)^2}{2}=0.
\end{equation}
in
$((s,\infty)\times\R)\setminus\omega$.
\end{lem}

Equivalently, $u(x,t)=\partial_x U(x,t)$  is an entropy  solution of
the Burgers equation~\eqref{eq:Burgers} with $f\equiv 0$ in
$((s,\infty)\times\R)\setminus\omega$. Of course, for each $t$,  $u(x,t)$ is a piecewise continuous
function of $x$, and at each of the countably many discontinuity points it makes a negative jump.
Since a velocity field determines its potential uniquely up to an additive constant we can also introduce
dynamics on velocity fields. We can introduce the space $\HH'$ of functions $w$ (actually, classes of equivalence of
functions since we do not distinguish two functions coinciding almost everywhere) such that
for some function $W\in\HH$ and almost every $x$,  $w(x)=W'(x)$.
This allows us to introduce the Burgers dynamics. We will
say that $ w_2=\Psi^{s,t}_\omega w_1$ if $w_1=W'_1$, $w_2=W'_2$, and $W_2=\Phi^{s,t}_\omega
W_1$  for some $W_1,W_2\in\HH$.

Often in the context of the Burgers dynamics, the functions in $\HH'$ will have negative
jump discontinuities, so called shocks. Although it is not essential, we can require the functions in $\HH'$ to be
right-continuous.

Let us denote $\HH(v_-,v_+)=\{W\in\HH:\ \lim_{x\to \pm\infty} (W(x)/x)=v_\pm \}$. The spaces $\hat\HH(v_-,v_+)$ are
defined
as classes of potentials in $\HH(v_-,v_+)$ coinciding up to an additive constant.

The following result shows that these spaces are invariant under HJHLO dynamics.
Along with Lemma~\ref{lem:cocycle} it allows to treat the dynamics as
a random dynamical system with perfect cocycle property (see, e.g.,~\cite[Section~1.1]{Arnold:MR1723992}).

\begin{lem} \label{lem:invariant_spaces}There is a set $\Omega\in\Fc_1$  with the following properties:
$\P_1(\Omega)=1$; for any $t\in\R$, $\theta^t\Omega=\Omega$;\quad  if $\omega\in\Omega$ then
for any $s,t$ with $s<t$,
\begin{enumerate}
 \item If $W\in\HH$, then $\Phi^{s,t}_\omega W\in\HH$.
 \item If $W\in\HH(v_-,v_+)$ for some  $v_-,v_+$, then $\Phi^{s,t}_\omega W\in\HH(v_-,v_+)$.
\end{enumerate}
\end{lem}

We denote the restrictions of $\Fc_1$ and $\P_1$ onto $\Omega$ by $\Fc$ and $\P$.
Since this restriction $\P$ still defines a Poisson point field with the same intensity measure, from
now on for convenience we remove from $\Omega_1$ the zero measure complement to $\Omega$ and work with probability
space $(\Omega,\Fc,\P)$.

In nonrandom setting the family of operators $(\Phi^{s,t})$ constructed via a variational
problem of type~\eqref{eq:optimization_problem}  is called a HJHLO evolution
semi-group, see~\cite[Definition 7.33]{Villani}, \cite{Fathi},
but in our setting it would be more precise to
call it a HJHLO cocycle.

Having defined the dynamics, we now turn to the main results.

\section{Main results}\label{sec:main_results}
 We say that $u(x,t)=u_\omega(x,t)$ is a global solution for the
cocycle $\Psi$ if there is a set $\Omega'$ with $\P(\Omega')=1$ such that for all
$\omega\in\Omega'$, all $s$ and $t$ with $s<t$, we have $\Psi^{s,t}_\omega u_\omega(\cdot,s)= u_\omega(\cdot,t)$.
We can also introduce the global solution as a skew-invariant function: $u_\omega(x)$ is called skew-invariant
if there is a set $\Omega'$ with $\P(\Omega')=1$ such that for any $t\in\R$, $\theta^t\Omega'=\Omega'$, and for any
$t>0$ and $\omega\in\Omega'$,
$\Psi^t_\omega u_\omega =u_{\theta^t\omega}$. If $u_\omega(x)$ is a skew-invariant function, then
$u_\omega(x,t)=u_{\theta^t\omega}(x)$ is a global solution. One can naturally view the potentials of $u_\omega(x)$ and
$u_\omega(x,s)$ as a  skew-invariant
function and global solution for the cocycle $\hat\Phi$.

Our first result is the description of global solutions.

\begin{thm}\label{thm:global_solutions}
For every $v\in\R$
there is a unique (up to zero-measure modifications) skew-invariant function
$u_v:\Omega\to\HH'$ such that for almost every $\omega\in\Omega$,
the potential $U_{v,\omega}$ defined by $U_{v,\omega}(x)=\int^x u_{v,\omega}(y)dy$ belongs to $\hat
\HH(v,v)$.

The potential $U_{v,\omega}$ is a unique skew-invariant potential in $\hat\HH(v,v)$. The skew-invariant functions
$U_{v,\omega}$ and $u_{v,\omega}$ are measurable  w.r.t.\ $\Fc|_{\R\times(-\infty,0]}$, i.e., they depend only on the
history of the forcing.  With probability 1, the realizations of $(u_{v,\omega}(y))_{y\in\R}$ are piecewise
linear with negative jumps between linear pieces. The spatial random process $(u_{v,\omega}(y))_{y\in\R}$ is stationary
and mixing.
\end{thm}
%
%
%
%
%
%

\begin{remark}\label{rem:stationary-solution}\rm Notice that this theorem can be interpreted as a 1F1S Principle: for
any velocity value $v$, the
solution at time $0$ with mean velocity $v$ is uniquely determined by
the history of the forcing: $u_{v,\omega}\stackrel{\rm a.s.}{=}\chi_v(\omega|_{\R\times (-\infty,0]})$ for some
deterministic functional $\chi_v$ of the point configurations on half-plane $\R\times (-\infty,0]$ (we actually describe
$\chi_v$ in the proof). Since the forcing is stationary in time, we obtain that $u_{v,\theta^t\omega}$ is a stationary
process in $t$, and the distribution of $u_{v,\omega}$ is an invariant distribution for the corresponding Markov
semi-group, concentrated on $\HH'(v,v)$.
\end{remark}
\medskip

The next result shows that each of the global solutions constructed in Theorem~\ref{thm:global_solutions} plays the
role of a one-point pullback attractor.
 To describe the domains of attraction  we will make assumptions on initial potentials
$W\in\HH$. Namely, we will assume that there is $v\in\R$ such that $W$ and $v$ satisfy one of the following sets of
conditions:
\begin{align}
v&=0,\notag\\
\liminf_{x\to+\infty} \frac{W(x)}{x}&\ge 0,  \label{eq:no_flux_from_infinity}\\
\limsup_{x\to-\infty} \frac{W(x)}{x}&\le 0,\notag
\end{align}
or
\begin{align}
v&> 0,\notag \\
\lim_{x\to-\infty} \frac{W(x)}{x}&= v,\label{eq:flux_from_the_left_wins}\\
\liminf_{x\to+\infty} \frac{W(x)}{x}&> -v,\notag
\end{align}
or
\begin{align}
v&< 0,\notag\\
\lim_{x\to+\infty} \frac{W(x)}{x}&= v,\label{eq:flux_from_the_right_wins}\\
\limsup_{x\to-\infty} \frac{W(x)}{x}&< -v.\notag
\end{align}

Condition~\eqref{eq:no_flux_from_infinity} means that there is no macroscopic flux of particles from infinity toward
the origin for the initial velocity profile $W'$. In particular, any $W\in\HH(0,0)$ or any $W\in\HH(v_-,v_+)$ with
$v_-\le 0$ and $v_+\ge 0$ satisfies~\eqref{eq:no_flux_from_infinity}. It is natural to call the arising phenomenon
a rarefaction fan.
We will see that in this case the long-term behavior is described
by the global solution $u_0$ with mean velocity $v=0$.

Condition~\eqref{eq:flux_from_the_left_wins} means that the initial velocity profile $W'$ creates the influx of
particles from $-\infty$ with effective velocity $v\ge 0$, and the influence of the particles at $+\infty$ is not as
strong.
In particular, any $W\in\HH(v,v_+)$ with
$v\ge 0$ and $v_+> -v$ (e.g., $v_+=v$) satisfies \eqref{eq:flux_from_the_left_wins}. We will see that in this case the
long-term
behavior is described by the global solution $u_v$.

Condition~\eqref{eq:flux_from_the_right_wins} describes a situation symmetric to~\eqref{eq:flux_from_the_left_wins},
where
in the long run the system is dominated by the flux  of particles from $+\infty$.

The following precise statement supplements Theorem~\ref{thm:global_solutions} and describes the basins of attraction of
the global solutions $u_v$ in terms of
conditions~\eqref{eq:no_flux_from_infinity}--\eqref{eq:flux_from_the_right_wins}.

\begin{thm}\label{thm:pullback_attraction} There is a set $\Omega''\in\Fc$ with $\P(\Omega'')=1$ such that if
$\omega\in\Omega''$, $W\in \HH$, and one of
conditions~\eqref{eq:no_flux_from_infinity},\eqref{eq:flux_from_the_left_wins},\eqref{eq:flux_from_the_right_wins}
holds:
then $w=W'$ belongs to the domain of pullback attraction of $u_v$ in the following sense:
for any $t\in\R$ and any $R>0$ there is $s_0=s_0(\omega)<t$ such that for all $s<s_0$
\[
 \Psi^{s,t}_\omega w(x) = u_{v,\omega}(x,t),\quad  x\in[-R,R].
\]
In particular,
\[
 \P\left\{\Psi^{s,t}_\omega w\big|_{[-R,R]}=u_{v,\omega}(\cdot,t)\big|_{[-R,R]}\right\}\to 1,\quad s\to-\infty.
\]

\end{thm}
\begin{remark}\rm The last statement of the theorem implies that for every $v\in\R$, the invariant
measure on $\HH'(v,v)$ described in Remark~\ref{rem:stationary-solution} is unique and for any initial condition
$w=W'\in\HH'$ satisfying one of
conditions~\eqref{eq:no_flux_from_infinity},\eqref{eq:flux_from_the_left_wins}, and~\eqref{eq:flux_from_the_right_wins},
the distribution of the random velocity profile at time $t$  converges to the unique stationary
distribution on $\HH'(v,v)$ as $t\to\infty$. However, our approach does not produce any convergence rate estimates.
\end{remark}
\begin{remark}\rm
Using space-time Galilean transformations, it is easy to obtain a version of
Theorem~\ref{thm:pullback_attraction} for attraction in a coordinate frame moving with constant velocity, but we omit
it for brevity.
\end{remark}
\medskip

The proofs of Theorems~\ref{thm:global_solutions} and~\ref{thm:pullback_attraction} are
given
in Sections~\ref{sec:global_solutions} and~\ref{sec:attractor}, but most of the preparatory work is carried out in
Sections~\ref{sec:subadd}, \ref{sec:concentration}, and \ref{sec:one-sided_minimizers}.

The long-term behavior of the cocycles $\Phi$ and $\Psi$ defined
through the optimization
problem~\eqref{eq:optimization_problem} depends on the asymptotic behavior of the action minimizers over long time
intervals.
The natural notion that plays a crucial role in this paper is the notion of backward one-sided infinite minimizers or
geodesics. A curve $\gamma:(-\infty,t]\to\R$ with $\gamma(t)=x$ is called a backward minimizer if its restriction onto
any time interval $[s,t]$ provides the minimum to the action $A_\omega^{s,t}(W,\cdot)$ defined in~\eqref{eq:action}
among paths connecting $\gamma(s)$ to $x$.

It can be shown (see Lemma~\ref{lem:geodesic direction}) that any backward minimizer $\gamma$ has an
asymptotic slope $v=\lim_{t\to
-\infty}(\gamma(t)/t)$. On the other hand, for every space-time point $(x,t)$ and every $v\in \R$ there is a backward
minimizer with slope $v$ and endpoint $(x,t)$. The following theorem describes the most important properties of
backward minimizers associated with the Poisson point field.

\begin{thm}
 For every $v\in\R$ there is a set of full measure $\Omega'$ such that for all $\omega\in\Omega'$ and any
$(x,t)\in\omega$ there is a unique
backward minimizer with asymptotic slope $v$. For any $(x_1,t_1),(x_2,t_2)\in\omega$ there is a time
$s\le \min\{t_1,t_2\}$ such that both minimizers coincide before $s$, i.e.,  $\gamma_1(r)=\gamma_2(r)$ for
$r\le s$.
\end{thm}

The proof of this core statement of this paper is spread over
Sections~\ref{sec:subadd} through~\ref{sec:one-sided_minimizers}.
In Section~\ref{sec:subadd} we apply the sub-additive ergodic theorem to derive the linear growth of action. In
Section~\ref{sec:concentration} we prove quantitative estimates on deviations from the linear growth. We use these
results in Section~\ref{sec:one-sided_minimizers} to analyze deviations of optimal paths from straight lines and deduce
the existence of infinite one-sided optimal paths and their
properties.

\section{Optimal action asymptotics and the shape function}\label{sec:subadd}

In this section we study the asymptotic behavior of the optimal action between space-time points $(x,s)$ and $(y,t)$
denoted by
\begin{align}\label{eq:optimal_action_between_two_points}
  A^{s,t}(x,y)=A^{s,t}_\omega(x,y)&=\min_{\gamma:\gamma(s)=x,\gamma(t)=y}(A^{s,t}_\omega(0,\gamma))
\\ &=
\min_{\gamma:\gamma(s)=x,\gamma(t)=y}\left(S^{s,t}(\gamma) - \omega^{s,t}(\gamma)\right)
\notag
\end{align}
(the minimum is taken over all absolutely continuous paths $\gamma$ or, equivalently, over all piecewise linear
paths with vertices at configuration points). Although to construct stationary solutions for the Burgers equation, we
will need the asymptotic behavior as $s\to{-\infty}$, it is more convenient and equally useful (due to the obvious
symmetry in the variational problem)
to
formulate most results for the limiting behavior as $t\to\infty$, and so we will here and in the next two sections.

We begin with some simple observations on Galilean shear transformations of the point field.
\begin{lem}\label{lem:shear} Let $a,v\in\R$ and let $L$ be a transformation of space-time defined by
$L(x,s)=(x+a+vs,s)$.
\begin{enumerate}
 \item
Suppose that  $\gamma$ is a path defined on a time interval
$[t_0,t_1]$ and let $\bar \gamma$ be defined by $(\bar\gamma(s),s)=L(\gamma(s),s)$.
Then
\[
 S^{t_0,t_1}(\bar\gamma)=S^{t_0,t_1}(\gamma)+(\gamma(t_1) -\gamma(t_0))v+\frac{(t_1-t_0)v^2}{2}.
\]
\item Let $L(\omega)$ be the point configuration obtained from $\omega\in\Omega$ by applying~$L$ pointwise. Then
$L(\omega)$ is also a Poisson process with Lebesgue intensity measure.

\item  Let $\omega\in\Omega$. For any time interval $[t_0,t_1]$ and any points $x_0,x_1,\bar x_0, \bar x_1$ satisfying
$L(x_0,t_0)=(\bar x_0, t_0)$ and $L(x_1,t_1)=(\bar x_1, t_1)$,
\[
A_{L(\omega)}^{t_0,t_1}(\bar x_0, \bar x_1)= A_{\omega}^{t_0,t_1}(x_0, x_1)+(x_1-x_0)v+\frac{(t_1-t_0)v^2}{2},
\]
and $L$ maps minimizers realizing $A_{\omega}^{t_0,t_1}(x_0, x_1)$ onto minimizers
realizing $A_{L(\omega)}^{t_0,t_1}(\bar x_0, \bar x_1)$.

\item For any points $x_0,x_1,\bar x_0,\bar x_1$ and any time interval $[t_0,t_1]$,
\[
  A^{t_0,t_1}(\bar x_0, \bar x_1)\stackrel{distr}{=} A^{t_0,t_1}(x_0,
x_1)+(x_1-x_0)v+\frac{(t_1-t_0)v^2}{2},
\]
where
\[
 v=\frac{(\bar x_1-x_1)-(\bar x_0-x_0)}{t_1-t_0}.
\]

\end{enumerate}
\end{lem}
\bpf The first part of the Lemma is a simple computation:
\begin{align*}
 S^{t_0,t_1}(\bar\gamma)& =\frac{1}{2}\int_{t_0}^{t_1}(\dot\gamma(s)+v)^2ds\\
 &=\frac{1}{2}\int_{t_0}^{t_1}\dot\gamma^2(s)ds+\int_{t_0}^{t_1}\dot\gamma(s)vds+\frac{1}{2}\int_{t_0}^{t_1}v^2ds.
\end{align*}
The second part holds since $L$ preserves the Lebesgue measure.
The third part follows from the first one since the images of paths transformed by $L$ are also paths passing through
the $L$-images of configuration points. The last part is a consequence of the previous two parts, since the
appropriate Galilean transformation sending $(x_0,t_0)$ to $(\bar x_0, t_0)$ and $(x_1,t_1)$ to $(\bar x_1, t_1)$
preserves the Lebesgue measure and the distribution of the Poisson process.
\epf

\medskip

The next useful property is the sub-additivity of action along any direction: for any velocity $v\in
\R$, and any $t,s\geq 0$, we have
\[ A^{0,t+s}(0,v(t+s))\leq A^{0,t}(0,vt) + A^{t,t+s}(vt,v(t+s)).\]
This means that we can apply Kingman's sub-additive ergodic theorem to the function $t\mapsto A^{0,t}(0,vt)$ if we can
show that $-\E A^{0,t}(0,vt)$ grows at most linearly in $t$. We claim this linear bound in
the following proposition:
\begin{lem}\label{lem:linbound}
Let $v\in \R$. There exist constants $C=C(v)>0$ and $t_0>0$ such that for all $t\geq t_0$
\[ \E |A^{0,t}(0,vt)| \leq Ct.\]
\end{lem}
\bpf Lemma~\ref{lem:shear} implies that it is enough to prove this for $v=0$. So in this proof we work with
$A^t=A^t(0,0)$.

Let $\gamma:[0,t]\to \R$ be a path realizing $A^t$. We have $\gamma(0)=\gamma(t)=0$. Let us split up
$\R^2$ into unit blocks
$B_{i,j}=[i,i+1)\times [j,j+1)$, for $i,j\in \ZZ$. We define $\calA$ as the union of all indices $(i,j)$ such that
$\gamma$ passes through $B_{i,j}$. The set $\calA$ is a lattice animal, i.e., it is a connected
set that
contains the origin $(0,0)\in\ZZ^2$ (see, e.g.,\cite{Gandolfi-Kesten:MR1258174}). Let us introduce the event $E_{n,t} =
\{ \#\calA = n\}$.

\begin{lem}\label{lem:latani} There
are constants $C_1, C_2, R, t_0>0$ such that if $t\geq t_0$ and $n\geq R t$, then
\[ \P(E_{n,t}) \leq C_1\exp(-C_2n^2/t).\]
\end{lem}

\bpf
We define $X_{i,j}=\omega(B_{i,j})$, the number of Poisson points in $B_{i,j}$.
Define the weight of the animal $\calA$ as
\[ N_\calA = \sum_{\nu\in \calA} X_\nu.\]
Clearly, the number of Poisson points picked up by $\gamma$ between $0$ and $t$, is
upper bounded by $N_\calA$. Define $k_j=\#\{i\in\ZZ\ :\ (i,j)\in\calA\}$, the number of blocks hit on the $j^{\rm th}$
row. These blocks will form a connected row of length $k_j$, and the kinetic action accumulated between
$j$ and $(j+1)\wedge t$ can therefore be bounded
by
\[ \frac12\int_j^{(j+1)\wedge t} \dot{\gamma}^2(s)\,ds \geq \frac12\left(k_j-2\right)_+^2.\]
Here, $a_+=\max(0,a)$. This leads to the following bound on the action:
\[ A^t \geq \frac12\sum_{0\le j<t}\left(k_j-2\right)_+^2 - N_\calA.\]
On $E_{n,t}$ we have $\sum_{0\le j<t} k_j = n$. Since $a\mapsto (a-2)_+^2$ is convex, we can use Jensen's inequality to see that
\[ \frac12\sum_{0\le j< t}\left(k_j-2\right)_+^2 \geq \frac12 \ceilt \left(\frac{n}{\ceilt} -
2\right)_+^2.\]
Therefore,
\begin{equation}
 A^t \geq \frac12 \ceilt \left(\frac{n}{\ceilt} - 2\right)_+^2 - N_\calA.
\label{eq:upper_bound_on_animal_weight}
\end{equation}
We also know that $A^t\leq 0$ since we can use the identical zero path on $[0,t]$.
Hence, on
$E_{n,t}$ we have
\[ N_\calA \geq \frac12 \ceilt \left(\frac{n}{\ceilt} - 2\right)_+^2.\]
Furthermore, if $N_n$ is the weight of the greedy animal of size $n$ (i.e., the animal of size $n$ with greatest
weight), then $N_n\ge N_\calA$, and
\[ E_{n,t} \subset \left\{N_n \geq \frac12 \ceilt \left(\frac{n}{\ceilt} - 2\right)_+^2\right\}.\]

Let us recall that the reasoning in~\cite{Cox-Gandolfi-Griffin-Kesten:MR1241039} after
equation (2.12)
implies that, due to standard large deviation estimates and the exponential growth of the number of lattice 
animals as a function of size $n$, there are constants $K_1, K_2, y_0>0$, such that if
\begin{equation}
\label{eq:cond_on_y}
y\ge y_0,
\end{equation}
 then
\begin{equation}\label{eq:latanibound}
\P\{N_n \geq yn\} \leq K_1\exp(-K_2ny).
\end{equation}
We now need to make sure that~\eqref{eq:cond_on_y} holds for
$y= \frac{1}{2n} \ceilt \left(\frac{n}{\ceilt} - 2\right)_+^2$.
If we require
$n\geq \max(4, 8 y_0) \ceilt$, then
\[ \frac{1}{2n} \ceilt \left(\frac{n}{\ceilt} - 2\right)_+^2
= \frac{1}{2n} \ceilt \left(\frac{n-2\ceilt}{\ceilt}\right)_+^2
\ge \frac{1}{2n} \ceilt \left(\frac{n-\frac{n}{2}}{\ceilt}\right)^2
\geq \frac{1}{8} \frac{n}{\ceilt}\ge y_0,\]
and the lemma follows from~\eqref{eq:latanibound}. \epf

\begin{remark}\rm We will choose the constant $R$ to be an integer, making it larger if needed.
\end{remark}

From~\eqref{eq:upper_bound_on_animal_weight} we already know that on $E_{n,t}$ we have
$0\geq A^t \geq -N_n$.
We wish to use this to estimate $\E |A^t|$, but we need an extension of~\eqref{eq:latanibound}.
\begin{lem}\label{lem:ENn}
For any $k\geq 1$, there is $c_k>0$ such that for all $n\geq 1$,
\[ \E N_n^k \leq c_kn^k.\]
\end{lem}
\bpf
Clearly,
\[ \E N_n^k = \sum_{i=0}^{\lfloor y_0n\rfloor} i^k\P\{N_n=i\} + \sum_{i=\lfloor y_0n\rfloor+1}^{\infty} i^k\P\{N_n=i\}. 
\]
We can bound the first term simply by
\[ \sum_{i=0}^{\lfloor y_0n\rfloor} i^k\P\{N_n=i\} \leq (y_0n)^k.\]
For the second term we can use \eqref{eq:latanibound}:
\[ \sum_{i=\lfloor y_0n\rfloor+1}^{\infty} i^k\P\{N_n=i\}\leq \sum_{i=\lfloor y_0n\rfloor+1}^{\infty}
K_1i^k\exp(-K_2i).\]
The right-hand side is bounded in $n$
and
the proof is complete. \epf

Lemma~\ref{lem:linbound} now follows from Lemmas~\ref{lem:latani}
and~\ref{lem:ENn}:
\begin{eqnarray*}
\E |A^t| & = & \sum_{n\le Rt} \E |A^t|\ONE_{E_{n,t}} + \sum_{n> Rt } \E |A^t|\ONE_{E_{n,t}}\\
& \leq & \E N_{[ Rt ]} +  \sum_{n> Rt } \E N_n\ONE_{E_{n,t}}\\
& \leq & Rc_1t + \sum_{n>Rt} \sqrt{\E N_n^2}\sqrt{\P(E_{n,t})}\\
& \leq & Rc_1t + \sqrt{c_2}\sum_{n>Rt}  \sqrt{C_1} n\exp(-C_2n^2/(2t))\\
& \leq & Ct,
\end{eqnarray*}
for $C$ big enough. \epf

In fact, we can use the last calculation to obtain the following generalization of Lemma~\ref{lem:linbound}  for
higher moments of $A^t$:

\begin{lem} \label{lem:moments_of_At} Let $k\in\NN$. Then there are constants $C(k),t_0(k)>0$ such that
\[ \E (|A^t|^k)\le C(k) t^k,\quad t\ge t_0(k).\]
\end{lem}

\medskip

Now a standard application of the sub-additive ergodic theorem shows that there exists a \emph{shape function} $\alpha(v)$ such that
\begin{equation}
\frac{A^{0,t}(0,vt)}{t}\to\alpha(v),\quad \text{a.s.\  and
in\ } L^1,\quad t\to\infty.
\label{eq:def_of_shape_function}
\end{equation}
Furthermore, $\alpha(0)<0$, since $A^t\leq 0$ and $\alpha(0)\leq \E(A^t)<0$. It turns out that the shape function $\alpha(v)$ is quadratic in $v$:
\begin{lem}\label{lem:shape-function} The shape function satisfies
\[
 \alpha(v)=\alpha(0)+\frac{v^2}{2},\quad v\in\R.
\]
\end{lem}
\bpf The Galilean shear map $(x,t)\mapsto(x+vt,t)$ transforms the paths connecting $(0,0)$ to
$(0,t)$ into paths connecting $(0,0)$ to $(vt,t)$.
 Lemma~\ref{lem:shear} implies that
under this map the optimal action over these paths is altered by a deterministic
correction $v^2t/2$. Since $\alpha$ is a constant almost surely we obtain
the
statement of the lemma.
\epf

We know now from~\eqref{eq:def_of_shape_function} that  $A^{0,t}(0,vt)\sim\alpha(v) t$ as $t\to\infty$ with
probability~1. However, this is not enough for our purposes since we need quantitative estimates on deviations of
$A^{0,t}(0,vt)$ from~$\alpha(v) t$. This is the material of the next section.


\section{Concentration inequality for optimal action}\label{sec:concentration}
 The goal of this section is to prove a concentration inequality for
$A^t(vt)=A^t(0,vt)=A^t_\omega(0,vt)=A^{0,t}_\omega(0,vt)$:

\begin{thm}\label{thm:concentration_around_alphat} There are positive
constants $c_0,c_1,c_2,c_3,c_4$ such that for any $v\in\R$, all
$t>c_0$, and
all
$u\in(c_3t^{1/2}\ln^2 t, c_4t^{3/2}\ln t]$,
\[
 \P\{|A^t(0,vt)-\alpha(v) t |>u\}\le c_1\exp\left\{-c_2\frac{u}{t^{1/2}\ln t}\right\}.
\]
\end{thm}
Due to the invariance under shear transformations (Lemmas~\ref{lem:shear} and~\ref{lem:shape-function}), it is
sufficient to prove this theorem for $v=0$.
We will first derive a similar inequality with $\alpha(0)t$ replaced by $\E A^t$, and then we will have to
estimate the corresponding approximation error. We recall that $A^t\leq 0$.

\begin{lem}\label{lem:concentration_around_mean} There are positive constants $b_0,b_1,b_2,b_3$ such that for all
$t>b_0$ and all
$u\in(0, b_3t^{3/2}\ln t]$,
\[
 \P\{|A^t-\E A^t|>u\}\le b_1\exp\left\{-b_2\frac{u}{t^{1/2}\ln t}\right\}.
\]
\end{lem}

The method of proof is derived from that for the generalized Hammersley's process in~\cite{CaPi}, but we have to take into
account that the optimal paths are allowed to travel  arbitrarily far within any bounded time interval in search for areas rich with configuration points. However, the situation where they decline too far from the kinetically most efficient path
is not typical. In the remaining part of this section we will often use the following lemma  showing that with
high probability the minimizer $\gamma$ connecting $(0,0)$ to $(0,t)$ stays within
distance $Rt$ from the origin, where $R$ was introduced in Lemma~\ref{lem:latani}.

\begin{lem}\label{lem:path_in_wide_rectangle_whp} There is a constant $C_3$ such that if $t\ge t_0$ and $u\ge Rt$ then
\[
\P\left\{\max_{s\in[0,t]}|\gamma(s)|>u\right\}\le C_3\exp(-C_2 u^2/t),
\]
where constants $C_2,R,t_0$ were introduced in Lemma~\ref{lem:latani}.
\end{lem}
\bpf
If $\max\{|\gamma(s)|:\ s\in[0,t] \}>u$, then the size of the lattice animal
$\calA$ traced by $\gamma$ is at least $u$.  Lemma~\ref{lem:latani} implies
\[
 \P\left\{\max_{s\in[0,t]}|\gamma(s)|>u\right\}\le\sum_{n\ge u} C_1\exp(-C_2n^2/t)\le C_3\exp(-C_2 u^2/t)
\]
for a constant $C_3$, since the  first term of the series is bounded by
$C_1\exp(-C_2u^2/t)$ and the ratio of two consecutive terms is bounded by $\exp(-C_2 R)$.
\epf

\medskip

Having Lemma~\ref{lem:path_in_wide_rectangle_whp} in mind, we define $\tilde A^t$ to be the optimal action over all
paths connecting $(0,0)$ to $(0,t)$ and staying within $[-Rt,Rt]$.

\begin{lem}\label{lem:probability_that_Kesten_action_worse} Let constants $t_0,R,C_2,C_3$ be defined in
Lemmas~\ref{lem:latani} and~\ref{lem:path_in_wide_rectangle_whp}. For any $t>t_0$,
\[
 \P\{A^t\ne \tilde A^t\}\le  C_3\exp(-R^2C_2 t).
\]
\end{lem}
\bpf It is sufficient to notice that
\[
 \P\{A^t\ne \tilde A^t\}\le\P\left\{\max_{s\in[0,t]}|\gamma(s)|>Rt\right\}
\]
and apply Lemma~\ref{lem:path_in_wide_rectangle_whp}.
\epf

\begin{lem}\label{lem:difference_of_rectangle_expectation_from_true} There is a constant $D_1$ such that for all
$t>t_0$,
\[
 0\le \E \tilde A^t-\E A^t\le -\E A^t\ONE_{\left\{\sup_{s\in[0,t]}|\gamma(s)|>Rt\right\}}\le D_1.
\]
\end{lem}
\bpf
The first two inequalities are obvious, since we have that $0\geq \tilde A^t\geq A^t$. For the last one, we have
\begin{align*}
-\E A^t\ONE_{ \left\{\sup_{s\in[0,t]}|\gamma(s)|>Rt\right\}}
&\le \sum_{n>Rt} \E( N_n\ONE_{E_{n,t}})
\\&\le \sum_{n>Rt} \sqrt{\E N^2_n}\sqrt{\P(E_{n,t})}
\\&\le  \sum_{n>Rt} \sqrt{c_2}n\sqrt{C_1}\exp(-C_2n^2/(2t)),
\end{align*}
where we used Lemmas~\ref{lem:latani} and~\ref{lem:ENn}. The statement follows since the last series is uniformly
convergent for $t>t_0$.
\epf

To obtain a concentration inequality for $\tilde A$, we will apply the following lemma by Kesten
\cite{Kesten:MR1221154}:

\begin{lem}
Let $(\Fc_k)_{0\le k\le N}$ be a filtration and let $(U_k)_{0\le k\le N}$ be a family of nonnegative random
variables
measurable with respect to $\Fc_N$. Let $(M_k)_{0\le k\le N}$ be a martingale with respect to $(\Fc_k)_{0\le k\le N}$.
Assume that for some constant $c>0$ the increments $\Delta_k=M_k-M_{k-1}$ satisfy
\[
 |\Delta_k|<c,\quad k=1,\ldots,N,
\]
and
\[
 \E(\Delta_k^2|\ \Fc_{k-1}) \le \E(U_k|\ \Fc_{k-1}).
\]
Assume further that for some positive constants $c_1,c_2$ and some $x_0\ge e^2c^2$ we have
\[
 \P\left\{\sum_{k=1}^N U_k>x\right\}\le c_1\exp(-c_2x),\quad x\ge x_0.
\]
Then
\[
 \P\{M_N-M_0\ge x\}\le c_3\left(1+c_1+\frac{c_1}{c_2
x_0}\right)\exp\left(-c_4\frac{x}{x_0^{1/2}+c_2^{-1/3}x^{1/3}}\right),\quad x>0,
\]
where  $c_3,c_4$  are universal positive constants that do not depend on $N,c,c_1,c_2,x_0$, nor on the
distribution of $(M_k)_{0\le k\le N}$ and $(U_k)_{0\le k\le N}$.
In particular,
\[
 \P\{M_N-M_0\ge x\}\le c_3\left(1+c_1+\frac{c_1}{c_2 x_0}\right)\exp\left(-c_4\frac{x}{2\sqrt{x_0}}\right),\quad x\le
c_2x_0^{3/2}.
\]
\end{lem}

To use this lemma in our framework, we must introduce an appropriate martingale.
For a given $t$ we consider the rectangle $Q(t)=[-Rt,Rt]\times [0,t]$ and
partition it  into $N=2Rt\cdot t=2Rt^2$ disjoint unit squares: $Q(t)=\bigcup_{k=1}^{N} B_k$. The order of
enumeration is not important. Here we assume that $t\in\NN$, but it is easy to adapt the reasoning
to the case of non-integer $t$.

We introduce a filtration $(\calF_k)_{0\le k\le N}$ in the following way. We set
$\calF_0=\{\emptyset,\Omega\}$ and
\[
 \calF_k=\sigma\left(\omega\bigr|_{\bigcup_{j=1}^k B_j} \right),\quad k=1,\ldots,N.
\]
We  introduce a martingale $(M_k,\calF_k)_{0\le k\le N}$ by
\[
 M_k=\E (\tilde A^t|\calF_k),\quad 0\le k\le N.
\]

We denote by $P_k$ the distribution of $\omega\bigr|_{B_k}$ on the sample space $\Omega_k$ of finite
point configurations in $B_k$. For $\omega,\sigma\in\prod_{k=1}^{N}\Omega_k$ we write
\[
 [\omega,\sigma]_k=(\omega_1,\ldots,\omega_k,\sigma_{k+1},\ldots,\sigma_{N})\in\prod_{k=1}^N\Omega_k.
\]
Then
\begin{align*}
 \Delta_k(\omega_1,\ldots,\omega_k)&:=M_k-M_{k-1}\\
 &= \int \tilde A^t_{[\omega,\sigma]_k}\prod_{j=k+1}^N dP_j(\sigma_j)-\int
\tilde A^t_{[\omega,\sigma]_{k-1}}\prod_{j=k}^N
dP_j(\sigma_j)\\
&=\int \left(\tilde A^t_{[\omega,\sigma]_k}- \tilde A^t_{[\omega,\sigma]_{k-1}}\right)\prod_{j=k}^N
dP_j(\sigma_j).
\end{align*}


\begin{lem}
\label{lm:estimating_integrand_in_Kesten_martingale_difference}
Let $I_k$ denote the indicator that the minimizer connecting $(0,0)$ to $(0,t)$ and staying in $[-Rt,Rt]$ passes
through a
Poissonian point in~$B_k$. Then
\[
|\tilde A^t_{[\omega,\sigma]_k}- \tilde A^t_{[\omega,\sigma]_{k-1}}| \le
\max\{I_k([\omega,\sigma]_k),I_k([\omega,\sigma]_{k-1})\}\max\{\omega(B_k),\sigma(B_k)\}.
\]
\end{lem}
\bpf
Suppose we delete the points of $\omega$ in $B_k$. When we then consider the minimizer for $[\omega,\sigma]_k$,
we decrease the number of Poissonian points contributing to the action by at most $\omega(B_k)$, and only decrease the kinetic action. Comparing the resulting path with
the minimizer for $[\omega,\sigma]_{k-1}$, we obtain
\[
  \tilde A^t_{[\omega,\sigma]_{k-1}}\le \tilde A^t_{[\omega,\sigma]_k}+\omega(B_k).
\]
Similarly, we get
\[
\tilde A^t_{[\omega,\sigma]_{k}}\le \tilde A^t_{[\omega,\sigma]_{k-1}} + \sigma(B_k).
\]
This shows that
\[ |\tilde A^t_{[\omega,\sigma]_k}- \tilde A^t_{[\omega,\sigma]_{k-1}}| \le \max\{\omega(B_k),\sigma(B_k)\}.\]
Now remark that if none of the two minimizers (for $[\omega,\sigma]_k$ and $[\omega,\sigma]_{k-1}$) passes through a Poissonian point
inside $B_k$, then $\tilde A^t_{[\omega,\sigma]_k}$ and $\tilde A^t_{[\omega,\sigma]_{k-1}}$ coincide.
This completes the proof.
\epf

The next step is to define a truncated Poissonian configuration $\bar \omega$ by erasing all Poissonian points of
$\omega$ in each
block $B_j$ with $\omega(B_j)> b\ln t$, where $b>0$ will be chosen later. The restrictions of $\bar \omega$ to $B_j$,
$j=1,\ldots,N$ are
jointly independent.
Lemma~\ref{lm:estimating_integrand_in_Kesten_martingale_difference} applies to truncated configurations as well
and we obtain
\[
| \tilde A^t_{[\bar\omega,\bar\sigma]_k}- \tilde A^t_{[\bar\omega,\bar\sigma]_{k-1}}| \le b \ln t
\max\{I_k([\bar\omega,\bar\sigma]_k),I_k([\bar\omega,\bar\sigma]_{k-1})\},
\]
where $\bar\sigma$ is obtained from $\sigma$ in the same way as $\bar \omega$ from $\omega$.
Therefore,
\[
 |\Delta_k(\bar\omega_1,\ldots,\bar\omega_k)|\le
 b\ln t \int \max\{I_k([\bar\omega,\bar\sigma]_k),I_k([\bar\omega,\bar\sigma]_{k-1})\} \prod_{j=k}^N
dP_j(\sigma_j)\le b\ln t.
\]

We must now estimate the increments of the martingale predictable characteristic. This estimate is a
straightforward analogue of Lemma~4.3 of~\cite{CaPi}.
\begin{lem}
Let $U_k=2(b\ln t)^2I_k$. Then, with probability $1$, $|U_k(\bar\omega)|\le 2(b\ln t)^2$ and
\[
 \E(\Delta_k^2(\bar\omega_1,\ldots,\bar\omega_k)|\calF_{k-1})\le \E(U_k(\bar\omega)|\calF_{k-1}).
\]
\end{lem}
\bpf
\begin{align*}
\E &(\Delta_k^2(\bar\omega_1,\ldots,\bar\omega_k)|\calF_{k-1})
=\int\left(\int\left(\tilde A^t_{[\bar\omega,\bar\sigma]_k}-
\tilde A^t_{[\bar\omega,\bar\sigma]_{k-1}}\right)\prod_{j=k}^N dP_j(\sigma_j)\right)^2 dP_k(\omega_k)\\
&\le\int\left(\int
\max\{I_k([\bar\omega,\bar\sigma]_k),I_k([\bar\omega,\bar\sigma]_{k-1})\}\cdot b\ln t
\prod_{j=k}^N dP_j(\sigma_j)\right)^2 dP_k(\omega_k)\\
&\le\int\int
\max\{I_k([\bar\omega,\bar\sigma]_k),I_k([\bar\omega,\bar\sigma]_{k-1})\}\cdot (b\ln t)^2
\prod_{j=k}^N dP_j(\sigma_j) dP_k(\omega_k)\\
&\le\int\int
(I_k([\bar\omega,\bar\sigma]_k)+I_k([\bar\omega,\bar\sigma]_{k-1}))\cdot (b\ln t)^2
\prod_{j=k}^N dP_j(\sigma_j) dP_k(\omega_k)\\
&= \E(U_k(\bar\omega)|\calF_{k-1}).
\end{align*}
\epf

We have
\begin{equation}
 \sum_{k=1}^N U_k(\bar\omega) = 2(b\ln t)^2\sum_{k=1}^{N} I_k(\bar\omega).
\end{equation}

Since
\[
 \sum_{k=1}^{N} I_k(\bar\omega)\le \# \calA (\bar\omega),
\]
we can write
\[
 \P\left\{\sum_{k=1}^N U_k(\bar\omega)>x\right\}\le \P\left\{\#\calA  (\bar\omega)>\frac{x}{2(b\ln t)^2}\right\}\le
\sum_{n>x/(2(b\ln t)^2) }  \P\{\bar\omega\in E_{n,t}\}.
\]

It is easy to see that Lemma~\ref{lem:latani} applies to $\bar\omega$ as well as to $\omega$, since its proof depends
only on the tail
estimate for the number of configuration points in each block.
We can conclude that
\begin{align*}
 \P\{\bar\omega\in E_{n,t}\} \leq C_1\exp(-C_2n^2/t),\quad n\ge Rt,\ t\ge t_0,
\end{align*}
where $C_1,C_2,R,t_0$ were introduced in Lemma~\ref{lem:latani}.

Combining the last two inequalities and choosing $x_0=2Rt(b\ln t)^2$, we can write for $x>x_0$
\begin{align*}
  \P\left\{\sum_{k=1}^N U_k(\bar\omega)>x\right\}
&\le C_1\sum_{n>x/(2(b\ln t)^2)}\exp(-C_2 n^2/t)
\\&\le C_4 \exp(-C_2
x^2/(4t(b\ln t)^4))\\
&\le  C_4 \exp(-C_2 x x_0/(4t(b\ln t)^4))
\\&\le C_4 \exp(-C_5 (b\ln t)^{-2}x).
\end{align*}

The above estimates on $\Delta_k(\bar\omega)$ and $U_k(\bar\omega)$ allow to apply Kesten's lemma with
$c=2b\ln t$, $c_1=C_4$, $c_2=C_5(b\ln t)^{-2}$, $x_0=2Rt(b\ln t)^2$ and obtain the following statement:

\begin{lem}\label{lem:outcome_of_Kestens_lemma}There are constants
$C_6,C_7,C_8,t_0>0$ such that for $t>t_0$ and $x\le C_8bt^{3/2}\ln t$,
\[
\P\{|\tilde A^t(\bar\omega)-\E \tilde A_t(\bar\omega)|>x\}\le C_6\exp\left(-C_7\frac{x}{b t^{1/2}\ln t}\right).
\]
\end{lem}

\begin{lem}\label{lem:deviation_of_truncation_from_original}
With probability 1,
\[
 \tilde A^t(\omega)\le \tilde A^t(\bar \omega).
\]
Also, we can choose $b$ and $t_0$ such that for all $t>t_0$ and $x>0$,
 \[
  \P\{\tilde A^t(\bar \omega) -\tilde A^t(\omega)>x\}\le 2e^{-x}.
 \]
\end{lem}
\bpf
The first statement of the lemma is obvious, and we have
\[
 0\le \tilde A^t(\bar \omega) -\tilde A^t(\omega)\le \sum_{k=1}^N \omega(B_k)\ONE_{\{\omega(B_k)>b\ln t\}}.
\]
By Markov's inequality and mutual independence of $\omega|_{B_k}$, $k=1,\ldots,N$,
\[
 \P\left\{\sum_{k=1}^N \omega(B_k)\ONE_{\{\omega(B_k)>b \ln t\}}>x\right\}\le
e^{-x}\left[\E e^{\omega(B_1)\ONE_{\{\omega(B_1)>b\ln t\}}}\right]^N.
\]
The lemma will follow from
\begin{equation}\label{eq:lim_factor_in_Markov_equals_1}
\lim_{t\to\infty}\left[\E e^{\omega(B_1)\ONE_{\{\omega(B_1)>b \ln t\}}}\right]^{2Rt^2}=1,
\end{equation}
which is implied by
\[
 \E e^{\omega(B_1)\ONE_{\{\omega(B_1)>b \ln t\}}}\le 1+\frac{\E  e^{2\omega(B_1)}}{e^{b\ln t}}\le 1+\frac{\E
e^{2\omega(B_1)}}{t^b},
\]
if we choose $b>2$.
\epf

The only missing part in the proof of Lemma~\ref{lem:concentration_around_mean} is the following corollary
of Lemma~\ref{lem:deviation_of_truncation_from_original}:
\begin{lem}\label{lem:expectation_of_truncation-expectation_over_Kesten_rectangle}
There is a constant $D_2$ such that for all
$t>t_0$,
\[
0\le \E \tilde A^t(\bar\omega)-\E \tilde A^t(\omega) <D_2.
\]
\end{lem}
\bpf[Proof of Lemma~\ref{lem:concentration_around_mean}]
Lemmas~\ref{lem:difference_of_rectangle_expectation_from_true}
and~\ref{lem:expectation_of_truncation-expectation_over_Kesten_rectangle}
 imply that
 for $u>D_1+D_2$
\begin{align*}
 \P\{|A^t(\omega)-\E A^t(\omega)|>u\}\le& \P\{|A^t(\omega)-\tilde A^t(\omega)|>(u-D_1-D_2)/3 \}
\\&+\P\{|\tilde
A^t(\omega)-\tilde A^t(\bar\omega)|>(u-D_1-D_2)/3\}
\\&+\P\{|\tilde
A^t(\bar\omega)-\E\tilde A^t(\bar\omega)|>(u-D_1-D_2)/3\}.
\end{align*}
The Lemma follows from the estimates of the three terms provided by
Lemmas~\ref{lem:probability_that_Kesten_action_worse},~\ref{lem:outcome_of_Kestens_lemma},
and
\ref{lem:deviation_of_truncation_from_original}.\epf

The following lemma quantifies how the growth of  $-\E A^t$ deviates from the linear one under argument doubling. We
will use this lemma to find an estimate on $\E A^t-\alpha(0)t$ which makes it possible to fill the gap between
Lemma~\ref{lem:concentration_around_mean} and Theorem~\ref{thm:concentration_around_alphat}.

\begin{lem}\label{lem:doubling}
There is a number $b_0>0$ such that for any $t>t_0$,
\[
 0\le 2\E A^t - \E A^{2t} \le b_0t^{1/2}\ln^2t.
\]
\end{lem}
\bpf The first inequality follows from $A^{0,2t}(0,0)\le A^{0,t}(0,0)+A^{t,2t}(0,0)$. Let us prove the second one.

Let $\gamma$ be the minimizer from $(0,0)$ to $(0,2t)$. Then
\[
 A^{2t}\ge \min_{|x|\le 2Rt} A^{0,t}(0,x)+\min_{|x|\le
2Rt}A^{t,2t}(x,0)+A^{2t}\ONE_{\left\{\max_{s\in[0,2t]}|\gamma(s)|>2Rt\right\}}.
\]
Therefore, by symmetry with respect to~$t$ and Lemma~\ref{lem:difference_of_rectangle_expectation_from_true},
\begin{equation}
\E A^{2t}\ge  2\E \min_{|x|\le 2Rt} A^{t}(0,x) - D_1.
\label{eq:doubling+constant}
\end{equation}

For $k\in I_t=\{-2Rt,\ldots,2Rt-2,2Rt-1\}$, we define a unit square $B_k=[k,k+1]\times[t-1,t]$.

Let now $\gamma$ be the minimizer from $(0,0)$ to $(x,t)$, with $x\in[k,k+1]$ for some
$k\in I_t$. Denote $t'=\sup\{s\le t: \gamma(s)\notin B_k\}$ and $x'=\gamma(t')$.

If $x'< k+1$, then by
reconnecting $(x',t')$ to $(k,t)$ we obtain
\[
 A^t(k)\le A^{t'}(x')+1/2\le A^t(x)+\omega(B_k)+1/2.
\]
If $x'= k+1$, then by
reconnecting $(x',t')$ to $(k+1,t)$ we obtain
\[
 A^t(k+1)\le A^{t'}(x')\le A^t(x)+\omega(B_k).
\]
Therefore,
\[
 A^t(x)\ge \min\{A^t(k),A^t(k+1)\}-\omega(B_k)-1/2,
\]
and~\eqref{eq:doubling+constant} implies
\[
 \E A^{2t}\ge  2 \E\min_{k\in I_t} A^t(k)-\E \max_{k\in I_t}\omega(B_k)-1/2-D_1.
\]
The second term grows logarithmically in $t$. Hence, for some constant $c>0$,
\[
 \E A^{2t}\ge  2 \E\min_{k\in I_t} A^t(k)-c(\ln t+1).
\]
Lemma~\ref{lem:shear} implies
\[
 \min_{x} \E A^t(x)=\E A^t(0).
\]
Therefore,
\begin{align}
\E A^{2t}&\ge  2 \E\min_{k\in I_t} A^t(k)-c(\ln t+1)  \notag
\\&\ge 2\min_{k\in I_t}\E A^t(k) - 2\E X_t - c(\ln t+1) \notag
\\&\ge 2\E A^t - 2\E X_t - c(\ln t+1), \label{eq:doubling1}
\end{align}
where
\[
 X_t=\max_{k\in I_t} \{(\E A^t(k)-A^t(k))_+\}.
\]
For a constant $r$ to be determined later, we introduce the event
\[
 E=\{X_t\le r(\ln^2 t)\sqrt{t}\}.
\]
Then
\[
 X_t\le r(\ln^2 t)\sqrt{t}\ONE_{E}+ X_t\ONE_{E^c}.
\]
Therefore,
\begin{equation}
 \E X_t\le r(\ln^2 t)\sqrt{t}+\sqrt{\E (X_{t})^2\P(E^c)}.
\label{eq:M_t}
\end{equation}
Let us estimate the second term. According to Lemma~\ref{lem:shear}, the random variables $A^t(k)-\E A^t(k)$, $k\in
I_t$ have the same distribution, so replacing the maximum in the definition of $X_t^2$ with summation we
obtain
\begin{equation}
\E X_{t}^2\le 4Rt \E (A^t-\E A^t)_+^2\le 4Rt \E (A^t)^2\le Ct^3,
\end{equation}
for some $C>0$ and all $t$ exceeding some $t_0$, where we used Lemma~\ref{lem:moments_of_At} in the last inequality.

Also, Lemma~\ref{lem:concentration_around_mean} shows that
\begin{align}
 \P(E^c)&\le \sum_{k\in I_t}\P\left\{A^t(k)-\E A^t(k) > r(\ln^2 t)\sqrt{t}\right\}\notag
\\ &\le 4R t b_1\exp\{-b_2r \ln t\}.\label{eq:probability_of_D_complement}
\end{align}
We can now finish the proof by choosing $r$ to be large enough and combining estimates
\eqref{eq:doubling1}--\eqref{eq:probability_of_D_complement}.
\epf

\medskip

With this lemma at hand we can now use the following statement (\cite[Lemma 4.2]{HoNe}):
\begin{lem}Suppose the functions $a: \R_+\to\R$ and $g:\R_+\to\R_+$  satisfy
the following conditions:
$a(t)/t\to \nu\in\R$ and  $g(t)/t \to 0$ as $t\to\infty$, $a(2t)\ge 2a(t)-g(t)$
and $\psi\equiv \limsup_{t\to\infty} g(2t) /g(t)< 2$. Then, for any $c > 1/(2-\psi)$,
and for all large $t$,
\[a(t) \leq \nu t + cg(t).\]
\end{lem}
Taking $a(t)=\E A^t$, $\nu=\alpha(0)$, $g(t)=b_0t^{1/2}\ln^2t$, $\psi=\sqrt{2}$, $c=2$,
we conclude that for $b'_0=2b_0$ and large $t$,
\[
0 \le \E A^t-\alpha(0) t\le b_0't^{1/2}\ln^2t,
\]
and Theorem~\ref{thm:concentration_around_alphat} follows from this estimate,
Lemma~\ref{lem:concentration_around_mean}, and the shear invariance established in Lemma~\ref{lem:shear}.\epf


\section{Existence and uniqueness of semi-infinite minimizers}\label{sec:one-sided_minimizers}

In this section we will study the properties of geodesics: these are paths $\gamma:[t_1,t_2]\to \R$ such that for
all $s\leq t \in [t_1,t_2]$, $\gamma|_{[s,t]}$ is the path that minimizes the action $A^{s,t}(\gamma(s),\gamma(t))$. We
will closely follow ideas by Howard and Newman in \cite{HoNe} and by W\"uthrich in \cite{Wu}, adapting them to our
specific situation, as was done in \cite{CaPi}.

\subsection{$\delta$-Straightness} The goal of this section is to estimate deviations of geodesics from straight lines.
 We will need the curvature of the shape function found in Section
\ref{sec:subadd}. Remember that
\[ \lim_{t\to \infty} \frac{A^{0,t}(0,vt)}{t} = \alpha(v) = \alpha(0) + \frac12 v^2.\]
Define $\alpha_0=\alpha(0)$. We will extend $\alpha$ to a function of $\R^2$:
\[ \alpha(p) := \alpha_0p_2 + \frac12
\frac{p_1^2}{p_2}=p_2\left(\alpha_0+\frac{1}{2}\left(\frac{p_1}{p_2}\right)^2\right)\]
(in this section we often denote the space and time coordinates of a space-time point $p\in\R^2$ by $p_1$ and
$p_2$, respectively).
This means that for all $p\in\R^2$ with $p_2>0$,
\[ \lim_{t\to \infty} \frac{A^{0,tp_2}(0,tp_1)}{t} \stackrel{\rm a.s.}= \alpha(p).\]
We need a convexity estimate of this function $\alpha$. Define for $p\in \R\times \R^+$ and for $L>0$
\[ \calC(p,L) := \left\{ q\in\R^2:\ q_2\in (p_2,2p_2]\ \mbox{and}\ \left|\frac{q_2}{p_2}\,p_1 - q_1\right| \leq
L\right\}.\]
So $\calC(p,L)$ is a parallelogram of width $2L$ along $[p,2p]$  (for any two points $p,q$ on the plane,
$[p,q]$ denotes the straight line segment connecting these two points). We need to consider the side-edges of this
parallelogram:
\[ \partial_S\calC(p,L) := \left\{ q\in\R^2:\ q_2\in (p_2,2p_2]\ \mbox{and}\ \left|\frac{q_2}{p_2}\,p_1 - q_1\right| =
L\right\}.\]
For a speed $v>0$ we define
\[{\rm Co}(v) = \{ p\in\R\times \R^+:\ |p_1|\leq p_2v\}.\]
The following lemma will play the role of Lemma 2.1 in \cite{Wu}.
\begin{lem}\label{lem:convest} Let $N,v>0$, $\delta\in(0,1)$. There are constants $M=M(N,v,\delta)>N$,
and $C=C(N,v,\delta)>0$
such that 
if points $p,p',q,q'$ satisfy $p\in\Co(v)$, $p_2>M$, $q\in\partial_S\calC(p,p_2^{1-\delta})$,
$p'_1=p_1$, $|p'_2-p_2|\le N$, $q'_1=q_1$, $|q'_2-q_2|\le N$, and
$p'_2<q'_2$, then
\begin{equation}\label{eq:convest}
\alpha(q'-p') + \alpha(p') - \alpha(q') \ge  C p_2^{1-2\delta}.
\end{equation}
\end{lem}
\bpf Using the definition of $\alpha$ and identities $p'_1=p_1$, $q'_1=q_1$, we get
\[
\alpha(q'-p') + \alpha(p')- \alpha(q') =
 \frac12\frac{(q_1-p_1)^2}{q'_2-p'_2} + \frac12\frac{{p_1}^2}{p'_2}- \frac12\frac{{q_1}^2}{q'_2}.
\]
So, $\alpha(q'-p') + \alpha(p')- \alpha(q')$ is a quadratic function in $q_1$. The minimum of this function
equals $0$, and is attained at $\tilde{q}_1$, the number defined by
\[ \frac{\tilde{q}_1-p_1}{q'_2-p'_2} = \frac{\tilde{q}_1}{q'_2}.\]
This means that $(\tilde{q}_1,q'_2)$ is a multiple of $p'$ and 
\begin{equation}
\label{eq:general-convexity}
\alpha(q'-p') + \alpha(p')- \alpha(q') =\frac12\frac{p'_2}{q'_2(q'_2-p'_2)}\,(q_1-\tilde{q}_1)^2.
\end{equation}
Now note that
\[ \tilde{q}_1 = \frac{q'_2}{p'_2}\,p_1 = \frac{q_2}{p_2}\,p_1 + \left(\frac{q'_2}{p'_2} - \frac{q_2}{p_2}\right)p_1.\]
We can estimate the second term by
\begin{align*}
 \left| \left(\frac{q'_2}{p'_2} - \frac{q_2}{p_2}\right)p_1\right| &= \frac{|p_2q'_2-p'_2q_2|}{p_2p'_2}\,|p_1| 
\le \frac{p_2|q'_2-q_2|+ |p_2-p'_2|q_2}{p_2p'_2}\,|p_1|
\\ &\le N\frac{|p_1|}{p'_2}+2N\frac{|p_1|}{ p'_2}\le c_1
\end{align*}
for some constant $c_1=c_1(M,N,v)$  if  $M$ is chosen to be sufficiently large and $p_2>M$, since 
$p'_2\ge p_2-N$ and $|p_1|\le v p_2$.

Since $|q_1-q_2p_1/p_2|=p_2^{1-\delta}$, we have that for large $M$ and $p_2>M$,
\begin{equation}
\label{eq:estimate-of-the-discrepancy}
|\tilde{q}_1-q_1| \ge  p_2^{1-\delta} -  \left| \left(\frac{q'_2}{p'_2} - \frac{q_2}{p_2}\right)p_1\right|
 \ge p_2^{1-\delta} - c_1  \ge c_2 p_2^{1-\delta},
\end{equation}
for some constant $c_2=c_2(M,N,v,\delta)$.

Also,
\begin{equation}
\label{eq:coefficient-in-front-of-quadratic}
\frac{p'_2}{q'_2(q'_2-p'_2)}\ge\frac{p_2-N}{(2p_2+N)(p_2+2N)}\ge \frac{c_3}{p_2}.
\end{equation}
for some $c_3=c_3(M,N)$. The lemma now follows from
\eqref{eq:general-convexity}, \eqref{eq:estimate-of-the-discrepancy}, \eqref{eq:coefficient-in-front-of-quadratic}.
\epf

This deterministic convexity lemma, together with the concentration bound of Section~\ref{sec:concentration}, will help
us to show that geodesics cannot make large deviations from a straight line. For $p\in \R^2$, we define
\begin{equation}
K(p,R):= \left\{q\in\R^2:\ |q_1-p_1|\leq R\ \mbox{and}\ |q_2-p_2|\leq R\right\}.
\label{eq:K-square}
\end{equation}
For $p,q\in \R^2$ satisfying $p_2<q_2$, we denote by $\gamma_{p,q}$ the optimal path from $p$ to $q$. For
$p,z\in \R^2$ with $0<p_2<z_2$, we define the event
\[ G(p,z) = \left\{ \exists \tilde{0}\in K((0,0),1), \tilde{z}\in K(z,1):\ \gamma_{\tilde{0},\tilde{z}}\cap K(p,1)
\neq \emptyset\right\}.\]
This says that a geodesic starting close to $(0,0)$ and ending near $z$, passes close to $p$. To bound the probability
of the event $G(p,z)$, we first have to control the action to and from points close to $p$. For $p,z\in\R^2$ with
$p_2<z_2$, we denote
\[ A(p,z)=A^{p_2,z_2}(p_1,z_1).\]
\begin{lem}\label{lem:boundaction}
Suppose $p,z\in \R^2$. Let $\tilde{p}\in K(p,1)$, $\tilde{z}\in K(z,1)$ with $\tilde{p}_2<\tilde{z}_2$. Define $\underline{p}=(p_1,p_2-2)$, $\bar{p}=(p_1,p_2+2)$ and similarly $\underline{z}$ and $\bar{z}$. Then
\[A(\tilde{p},\tilde{z})\geq A(\underline{p},\bar{z}) - 1.\]
If $p_2+2<z_2-2$, then
\[ A(\bar{p},\underline{z})+1 \geq A(\tilde{p},\tilde{z}).\]
\end{lem}
\bpf
Let $\tilde{\gamma}$ be the optimal path (i.e., the path that picks up the least action) from  $\tilde{p}$ to
$\tilde{z}$. Define $\gamma$ as the path that starts at $\underline{p}$, moves at constant speed to
$\tilde{p}$, then follows $\tilde{\gamma}$,
and then moves at constant speed to $\bar{z}$. Then, denoting $A[\gamma]$ for the action picked up by a path $\gamma$,
we get
\begin{eqnarray*}
A(\underline{p},\bar{z}) & \leq & A[\gamma]\\
 & \leq & \frac12 + A(\tilde{p},\tilde{z}) + \frac12.
\end{eqnarray*}
We use that the first part of this path picks up at most $v^2/2$ action, where the speed $v\leq 1$. For the path from $\tilde{z}$ to $\bar{z}$ we get the same upper bound.

For the second inequality, note that
\[ A(\tilde{p},\tilde{z}) \leq A(\tilde{p},\bar{p}) + A(\bar{p},\underline{z}) + A(\underline{z},\tilde{z}).\]
Clearly, as described above, we have $A(\tilde{p},\bar{p})\leq 1/2$ by taking a similar path from $\tilde{p}$ to $\bar{p}$. Likewise, $A(\underline{z},\tilde{z})\leq 1/2$. This proves the lemma.
\epf

\begin{lem}\label{lem:dstraight1}
Fix $\delta\in(0,1/4)$ and $v>0$. There exist constants $c_1, c_2, M>0$, such that for all $p\in {\rm Co}(v)$ with $p_2>M$ and $z\in \partial_S\calC(p,p_2^{1-\delta})$, we have
\[ \P(G(p,z))\leq c_1\exp\left(-c_2p_2^{1/2-2\delta}/\log(p_2)\right).\]
\end{lem}
\bpf
Let us choose $M$ according to Lemma~\ref{lem:convest} and make sure that $M>8$ (we may always increase~$M$
 keeping Lemma~\ref{lem:convest} true). Take $z\in \partial_S\calC(p,p_2^{1-\delta})$. Define $\underline{p}=(p_1,p_2-2)$,
$\bar{p}=(p_1,p_2+2)$,
and likewise $\underline{0}$, $\bar{0}$, $\underline{z}$ and $\bar{z}$. The event $G(p,z)$ implies that there exist
three points  $\tilde{0}\in K(0,1)$, $\tilde{p}\in K(p,1)$ and $\tilde{z}\in K(z,1)$ with $\tilde{p}_2<\tilde{z}_2$ such
that
\[ A(\tilde{0},\tilde{z}) = A(\tilde{0},\tilde{p}) + A(\tilde{p},\tilde{z}).\]
It follows from Lemma \ref{lem:boundaction} that
\[ A(\bar{0},\underline{z})\geq A(\tilde{0},\tilde{z}) - 1,\ \  A(\underline{0},\bar{p})\leq A(\tilde{0},\tilde{p}) + 1\ \mbox{and}\ A(\underline{p},\bar{z})\leq A(\tilde{p},\tilde{z}) + 1.\]
Therefore,
\begin{equation}
 A(\bar{0},\underline{z})\geq A(\underline{0},\bar{p}) + A(\underline{p},\bar{z}) - 3.
 \label{eq:convexity-wrong-direction}
\end{equation}
Furthermore, consider $\alpha(\underline{z}-\bar{0})$. Since $p_2^{1-\delta}<p_2$, we have that $\bar z\in {\rm Co}(v+1)$.
\begin{eqnarray*}
\alpha(\underline{z}-\bar{0}) & = & \alpha_0(z_2-4) + \frac12\,\frac{z_1^2}{z_2-4}\\
& = & \alpha(\bar z) - 6\alpha_0 + \frac{3z_1^2}{(z_2+2)(z_2-4)}\\
& \leq & \alpha(\bar z) -6\alpha_0  + 6(v+1)^2.
\end{eqnarray*}
Here we use that $8<M\leq z_2$ and $|z_1|\leq (v+1)z_2$. Therefore, following the same reasoning, we can choose $L>0$ independent of $p\in {\rm Co}(v)$ such that
\begin{equation}
\label{eq:approximations-for-alpha}
\alpha(\underline{z}-\bar{0})\leq \alpha(\bar z) + L\quad \mbox{and}\quad \alpha(\bar{p}-\underline{0})\geq \alpha(\underline{p}) - L.
\end{equation}
Also, Lemma~\ref{lem:convest} implies
\[
\alpha(\bar{z}-\underline{p}) + \alpha(\underline{p}) \geq \alpha(\bar{z}) + C p_2^{1-2\delta}.
\]
Combining this with~\eqref{eq:convexity-wrong-direction} and~\eqref{eq:approximations-for-alpha}, we obtain
\[ ( A(\bar{0},\underline{z}) - \alpha(\underline{z}-\bar{0})) 
- (A(\underline{0},\bar{p}) - \alpha(\bar{p}-\underline{0})) - (A(\underline{p},\bar{z}) - \alpha(\bar{z}-\underline{p})) \geq -3  - 2L + Cp_2^{1-2\delta}.\]
Define the event corresponding to the third term on the l.h.s.
\[ E_3 = \left\{ -(A(\underline{p},\bar{z}) - \alpha(\bar{z}-\underline{p})) \ge \frac{C}{4}\,p_2^{1-2\delta}\right\},\]
and likewise $E_1$ and $E_2$. By enlarging $M$, we can make sure that $G(p,z)$ implies at least one of these three events. We will bound the probability of $E_3$, which is slightly more complicated than the other two, since $\bar{z}_2-\underline{p}_2$ cannot be made arbitrarily large by increasing $M$.

Using the shear transformation we know that
\[A(\underline{p},\bar{z}) - \alpha(\bar{z}-\underline{p})\stackrel{distr}{=} 
A((0,0),(0,\bar{z}_2-\underline{p}_2) - \alpha_0\cdot(\bar{z}_2-\underline{p}_2) =
A^{\bar{z}_2-\underline{p}_2}- \alpha_0\cdot(\bar{z}_2-\underline{p}_2).\]
Since $A^{\bar{z}_2-\underline{p}_2}
\ge A^{0,\bar{z}_2}(0,0)-A^{\bar{z}_2-\underline{p}_2,\bar{z}_2}(0,0)$ 
and $A^{\bar{z}_2-\underline{p}_2,\bar{z}_2}(0,0)\stackrel{distr}{=}A^{\underline{p}_2}$,
we obtain
\begin{align*}
\P(E_3)&\le 
\P\left\{-(A^{0,\bar{z}_2}(0,0) - \alpha_0 \bar{z}_2)+  (A^{\bar{z}_2-\underline{p}_2,\bar{z}_2}(0,0) - \alpha_0 \underline{p}_2 )\geq \frac{C}{4}\,p_2^{1-2\delta}\right\}
\\
&\leq \P\left\{|A^{\bar{z}_2} - \alpha_0 \bar{z}_2|\geq \frac{C}{8}\,p_2^{1-2\delta}\right\}
+  \P\left\{|A^{\underline{p}_2} - \alpha_0 \underline{p}_2 |\geq \frac{C}{8}\,p_2^{1-2\delta}\right\}.
\end{align*}
Theorem \ref{thm:concentration_around_alphat} guarantees the existence of constants $C_1,C_2>0$ such that for $M$ big enough and all $p$ with $p_2\geq M$,
\[ \P(E_3) \leq C_1\exp\left\{-C_2 |p_2|^{1/2-2\delta}/\log(p_2)\right\}.\]
Similar bounds hold for $\P(E_1)$ and $\P(E_2)$, proving the lemma.
\epf

The above Lemma can be used to show that a minimal path starting close to the origin and passing  close to $p$, with
high probability will not exit the slanted cylinder $\calC(p,p_2^{1-\delta})$ through the sides. Define the event
\[ G(p) = \left\{ \exists\ \tilde{0}\in K((0,0),1)\ \exists\ z\in \partial_S\calC(p,p_2^{1-\delta}):\
\gamma_{\tilde{0},z}\cap K(p,1) \neq \emptyset\right\}.\]

\begin{lem}\label{lem:dstraight2}
Fix $\delta\in(0,1/4)$ and $v>0$. There exist constants $c_1, c_2, \kappa, M>0$, such that for all $p\in {\rm Co}(v)$ with $p_2>M$ we have
\[ \P(G(p))\leq c_1\exp(-c_2p_2^\kappa).\]
\end{lem}
\bpf
Suppose $p\in\R^2$ with $p_2>M$. There exists a constant $c>0$ (depending on $v$) and points $z_1,\ldots, z_L \in \partial_S\calC(p,p_2^{1-\delta})$ with $L\leq cp_2$, such that
\[ \partial_S\calC(p,p_2^{1-\delta}) \subset \bigcup_{i=1}^L K(z_i,1).\]
This implies that
\[ G(p)\subset \bigcup_{i=1}^L G(p,z_i).\]
Therefore, by choosing $M$ large enough, $\kappa < 1/2-2\delta$ and using Lemma \ref{lem:dstraight1}, there exist $C_1, C_2, c_1, c_2>0$ such that
\[ \P(G(p)) \leq LC_1\exp(C_2p_2^{1/2-2\delta}/\log(p_2)) \leq c_1\exp(-c_2p_2^\kappa).\]
\epf

Now we are ready to prove $\delta$-straightness of geodesics, as was introduced by Newman in \cite{Ne}. For a path $\gamma$ and $t\in \R$, we define
\[ \gamma^{\rm out}(t) = \{ (\gamma(s),s)\ :\ s\geq t\},\]
which is the set of all points in the path $\gamma$ (more precisely, the \emph{graph} of $\gamma$), that are reached after time $t$. We also consider the following cone for all $x\in \R\times \R^+$ and $\eta>0$:
\begin{equation}
\label{eq:Co}
{\rm Co}(x,\eta) = \{ z\in \R\times \R^+\ :\ |z_1/z_2 - x_1/x_2| \leq \eta \},
\end{equation}
which is the cone starting in the origin of all points $z$ that have a corresponding speed closer than $\eta$ to the speed of $x$.

\begin{lem}[{\bf $\delta$-straightness}]\label{lem:delta_straightness}
For
$\delta\in(0,1/4)$ and $v>0$ we have with probability one that there exists $M>0$ (depending on $v$
and $\delta$) and $R>0$ (depending only on $\delta$), such that for all $\tilde{0}\in K((0,0),1)$, for all $z\in
\R\times \R^+$ and for all $p\in \gamma(\tilde{0},z)\cap {\rm Co}(v)$ with $p_2>M$, we have
\[ \gamma^{\rm out}(p_2) \subset {\rm Co}(p,Rp_2^{-\delta}),\]
for $\gamma=\gamma_{\tilde{0},z}$.
\end{lem}
This lemma states that if a geodesic starting near $(0,0)$ passes through a remote point $p$, it has to stay in a
narrow cone around the ray $\R^+\cdot p$.


\bpf
Consider the events $G(\bar{p})$ for all $\bar{p}\in \Z\times \Z^+\cap {\rm Co}(v')$, with $v'>v$. Using
Lemma~\ref{lem:dstraight2} and Borel--Cantelli Lemma, we can choose $M$ big enough, such that for all
$\bar{p}\in \Z\times \Z^+\cap {\rm Co}(v')$ with $\bar{p}_2\geq M$ the event $G(\bar{p})$ does not happen. Increase $M$
if necessary to ensure that if $p_2\geq M$,
\[ K(p,1+(p_2+1)^{1-\delta})\subset {\rm Co}(p,2p_2^{-\delta}).\]
So for any $p\in\R\times\R^+$ with $p_2\geq M$, we now know that if $\bar{p}\in K(p,1)$, then
\[ \calC(\bar{p},\bar{p}_2^{1-\delta})\subset {\rm Co}(p,2p_2^{-\delta}).\]
Let $\bar{0}\in K((0,0),1)$. Now suppose there exists $z\in\R\times\R^+$ and $p\in \gamma_{\bar{0},z}\cap {\rm Co}(v)$ with $p_2>M+1$. Define $\bar{p}=(\lfloor p_1\rfloor,\lfloor p_2\rfloor)$. Suppose $z$ lies outside of the slanted tube ${\cal C}(\bar{p},\bar{p}_2^{1-\delta})$. We know that $G(\bar{p})$ does not happen, and since $p\in K(\bar{p},1)$, this implies that we can define $p^{(1)}$ as the crossing of $\gamma^{\rm out}(p_2)$ with the top edge of the tube ${\cal C}(\bar{p},\bar{p}_2^{1-\delta})$, and that $\gamma_{p,p^{(1)}}$ lies inside this tube. Define $\bar{p}^{(1)}=(\lfloor p^{(1)}_1\rfloor,\lfloor p^{(1)}_2\rfloor)$. We can proceed in a similar way to construct $p^{(2)},\bar{p}^{(2)},p^{(3)},\ldots,\bar{p}^{(m)}$, where we finish whenever $z\in {\cal C}(\bar{p}^{(m)},(\bar{p}^{(m)}_2)^{1-\delta})$. We have to check that for all $k\leq m$, $\bar{p}^{(k)}\in {\rm Co}(v')$, but this will follow from the following considerations.

Note that for $1\leq k\leq m$, where we define $p^{(0)}=p$ and $\bar{p}^{(0)}=\bar{p}$,
\[ {\cal C}(\bar{p}^{(k)},(\bar{p}^{(k)}_2)^{1-\delta})\subset {\rm Co}(p^{(k)},2(p_2^{(k)})^{-\delta})\ \ {\rm and}\ \ p_2^{(k)}\geq 2(p_2^{(k-1)}-1)\geq \frac32 p_2^{(k-1)}.\]
This implies that the average speed of any vector in $y\in \gamma^{\rm out}(p_2)$ satisfies
\[ |y_1/y_2 - p_1/p_2| \leq \sum_{k=0}^m 2(p_2^{(k)})^{-\delta} \leq \sum_{k=0}^m 2\left(\frac32\right)^{-\delta
k}p_2^{-\delta}\leq Rp_2^{-\delta},\]
if we choose $R>0$ large enough (depending only on $\delta$). This also shows that we have to choose $v'>v+RM^{-\delta}$.
\epf

\begin{coro}\label{cor:delta-straightness-prob}
For $\delta\in(0,1/4)$ and $v>0$ there exists $M,R,\kappa,C_1,C_2>0$, such that when we define the event
\begin{eqnarray*}
G_n & = & \Big\{\exists \tilde{0}\in K((0,0),1), z\in \R\times \R^+, p\in\gamma(\tilde{0},z)\ \mbox{with}\ p_2>n\
\mbox{and}\ p\in {\rm Co}(v)\ :\\
& & \ \  \gamma^{\rm out}(p_2) \not\subset {\rm Co}(p,R p_2^{-\delta})\Big\},
\end{eqnarray*}
we have for $n\geq M$
\[ \P(G_n)\leq C_1e^{-C_2n^\kappa}.\]
\end{coro}
\bpf It follows directly from the proof of Lemma~\ref{lem:delta_straightness} that the event $G_{n}$ is a subset
of the event that there exists $\bar{p}^{(n)}\in \Z\times \Z^+\cap {\rm Co}(v')$ with $\bar{p}^{(n)}_2\geq n-1$ such
that the event $G(\bar{p}^{(n)})$ does happen. Here we choose $v'>v+RM^{-\delta}$. The probability of this event is
clearly bounded by $C_1e^{-C_2n^\kappa}$, for an appropriate choice of constants, simply by Lemma~\ref{lem:dstraight2}.
\epf

\subsection{Existence and uniqueness of semi-infinite minimizers}
With $\delta$-straightness in hand, we can prove some important properties of minimizing paths.
A semi-infinite minimizer starting at  $(x,t)\in\R^2$ is a path $\gamma:[t,\infty)\to \R$ such that
$\gamma(t)=x$ and the restriction of $\gamma$ to any finite time interval is
a minimizer. We call  $(x,t)$ the endpoint of $\gamma$.

\begin{lem}\label{lem:geodesic direction}
With probability one, all semi-infinite minimizers have an asymptotic slope (velocity, direction): for every
minimizer $\gamma$ there
exists $v\in \R\cup\{\pm \infty\}$  depending on $\gamma$ such that
\[ \lim_{t\to\infty} \frac{\gamma(t)}{t} = v.\]
\end{lem}
\bpf Let us fix a sequence $v_n\to \infty$. Using the translation invariance of the Poisson point field,
with probability one, for any $q\in \Z^2$ we can choose  a
corresponding sequence of constants $M_n(q)>0$ such that the statement in Lemma~\ref{lem:delta_straightness} holds for
the entire sequence, for paths starting in $K(q,1)$.

Let us take some one-sided minimizer $\gamma$. If $\gamma(t)/t\to+\infty$ or $-\infty$, then the desired statement
is automatically true. In the opposite case we have
\[ \liminf_{t\to\infty} \frac{|\gamma(t)|}{t} < \infty.\]
This implies that there exist $n\geq 1$ and a sequence $t_m\to \infty$ such that $|\gamma(t_m)|/t_m \leq v_n$. We define $y_m=(\gamma(t_m),t_m)$ and choose $q\in\Z^2$ such
that $y_1\in K(q,1)$. For $m$ large enough, we will have that $t_m>M_n(q)$ and, therefore,
\[ \gamma^{\rm out}(y_m)\subset  q + {\rm Co}(y_m-q,R|y_m-q|^{-\delta}),\]
for some constant $R>0$ and $m$ large enough. Clearly this implies that $\gamma$ must have a finite asymptotic
slope.
\epf

\begin{lem}\label{lem:exist-geodesic}
With probability one, for every $v\in\R$ and for every sequence $(y_n,t_n)\in \R^2$ with $t_n\to \infty$ and
\[ \lim_{n\to \infty} \frac{y_n}{t_n} = v,\]
and for every $x\in \R^2$, there exists a subsequence $(n_k)$ such that the minimizing paths
$\gamma_{x,(y_{n_k},t_{n_k})}$ are an increasing collection of paths that converge to a semi-infinite minimizer starting
at $x$ and with asymptotic speed equal to $v$.
\end{lem}
\bpf Without loss of generality, we can assume that $x\in K((0,0),1)$. We take a sequence $v_m\to \infty$ and then
choose $M_m\to \infty$ and $R>0$ such that the statement of Lemma~\ref{lem:delta_straightness} holds for every
triplet $(v_m,M_m,R)$. From these triplets we choose a triplet $(v_0,M,R)$ with $|v|<v_0-2RM^{-\delta}$ (note that $R$ only depends on $\delta$).

By going to a subsequence, we make sure that for all
$n\geq 1$, $t_n\geq M$, $t_n\uparrow \infty$ and for all $k>n$ we have $(y_k,t_k)\in {\rm
Co}((v,1),Rt_{n}^{-\delta})$. Consider the paths $\gamma_n = \gamma_{x,(y_n,t_n)}$. We claim that for each $n\geq 1$
and $k>n$, the path $\gamma_k$  lies in the cone ${\rm Co}((v,1),2Rt_n^{-\delta})$ for times larger than
$t_n$. In fact, if $\gamma_k$ visits a point $p$ outside this cone (but inside ${\rm Co}(v_0)$)
at time $p_2\geq t_n$, then $\gamma_k$ violates the $\delta$-straightness condition (the relevant cone through
$p$ will not intersect ${\rm Co}((v,1),Rt_n^{-\delta})$, and therefore it does not contain $y_k$). In particular,
this means that there exists a $C>0$ (independent of $n$) such that all paths $\gamma_k$ with $k>n$ cross the
segment
\[J_n=[vt_n-Ct_n^{1-\delta},vt_n+Ct_n^{1-\delta}]\times \{t_n\}.\]
We claim that there exists a point $x_n\in J_n$ visited by an infinite number of paths $\gamma_k$.
With this claim in hand, we first choose $x_1$ and a subsequence
$k_1(n)$ such that every $\gamma_{k_1(n)}$ passes through $x_1$, then we choose $x_2$ in the segment
$$[vt_{k_1(1)}-Ct_{k_1(1)}^{1-\delta},vt_{k_1(1)}+Ct_{k_1(1)}^{1-\delta}]\times \{t_{k_1(1)}\}$$
such that an infinite number (a subsequence $k_2(n)\subset k_1(n)$)
 of the paths
$\gamma_{k_1(n)}$ pass through $x_2$,  and so on. The paths $\{\gamma_{x,x_n}\ :\ n\geq 1\}$ are then an
increasing collection of paths, and their union will be a semi-infinite minimizer starting at $x$, with asymptotic speed
$v$.\\

What remains is to show that an infinite number of the paths $\gamma_k, k>n$, pass through the same point of  the
segment
$$[vt_n-Ct_n^{1-\delta},vt_n+Ct_n^{1-\delta}]\times \{t_n\}.$$
 We already know that all these paths remain in the cone ${\cal C} = {\rm
Co}((v,1),Rt_n^{-\delta})$. We use the following fact about the Poisson process, which is easily checked using the
Borel--Cantelli Lemma: with probability one, for all $K>0$ there exists $T_0>0$ depending on $K$, such that for all
$T\geq T_0$ and all $x\in [-KT,KT]\times [-1,2T]$, there is at least one Poisson point in the ball of radius $T^{1/4}$
around $x$. Let us choose $K>0$ and $T>t_n$ such that ${\cal C}\cap (\R\times [0,t_n+T])\subset [-KT,KT]\times [-1,2T]$. Then we
choose $T_0>t_n$ according to this $K$ and take $T\geq T_0$.

Consider $k$ such that $t_k>t_n+2T$. Suppose
$\gamma_k$ did not pick up any Poisson point in the time interval $[t_n,t_n+2T]$. That implies that it would have some
constant speed $u$ in that interval. We also know that within distance $T^{1/4}$ of $(\gamma_k(t_n+T),t_n+T)$ there
will be some Poisson point, since $(\gamma_k(t_n+T),t_n+T)\in {\cal C}$. Define the path $\tilde{\gamma}_k$ that picks
up this Poisson point by moving to this point with constant velocity from the point $(\gamma_k(t_n),t_n)$ and then
moving with constant speed to $(\gamma_k(t_n+2T),t_n+2T)$; the rest of the time $\tilde{\gamma}_k$ coincides with
$\gamma_k$. When we consider the difference in action picked up by the two paths, we note that they start and end at the
same point in the time interval $[t_n,t_n+2T]$. This means that if we define $\delta(t) = \dot{\tilde{\gamma}}_k(t) -
\dot\gamma_k(t)$, then $\int_{t_n}^{t_n+2T} \delta(s)\,ds = 0$. Furthermore, there exists a constant $d>0$ depending
only on $v$ and the cone $\cal C$ such that $|\delta(s)|\leq dT^{-3/4}$. If we choose $T>d^4$, this leads to
\begin{eqnarray*}
A[\tilde{\gamma}_k] - A[\gamma_k] & = & -1 + \frac12 \int_{t_n}^{t_n+2T} ((u+\delta(s))^2 - u^2)\,ds\\
& = & -1 + \frac12 \int _{t_n}^{t_n+2T} \delta(s)^2\,ds\\
& \geq & -1 + d^2 T^{-1/2}\\
& < & 0.
\end{eqnarray*}
This contradicts the optimality of $\gamma_k$, and we conclude that  $\gamma_k$ must pick a
Poisson point in the time interval $[t_n,t_n+2T]$. The number of paths $\gamma_k$ is infinite, and there are only
finitely many Poisson points in the set ${\cal C}\cap \R\times [t_n,t_n+2T]$. Therefore, at time
$t_n$, an infinite number of paths $\gamma_k$ cross $J_n$ at the same point. Here we
use that when two minimizers meet at two Poisson points at distinct times, they actually will coincide for all times between these two
times (since minimizing paths between two Poisson points are almost surely unique).
\epf

\begin{lem}\label{lem:geodesics-do-not-share-more-than-one} With probability 1, the following statement holds: if
$\gamma_1$ and $\gamma_2$ are two (finite-time) geodesics, starting at the same Poisson point $p$, and for some $t>p_2$
we have $\gamma_1(t)<\gamma_2(t)$, then for all (relevant) $s>t$ we have $\gamma_1(s)<\gamma_2(s)$.
\end{lem}
\bpf The probability that there are two Poisson point connected by two distinct geodesics is zero. So we only have to
consider the situation where two paths with vertices $p,p_1,\ldots,p_n$ and $p,q_1,\ldots,q_m$ intersect
transversally,
i.e., for some $k,j$, $[p_j,p_{j+1}]\cap [q_k,q_{k+1}]=\{x\}$ for some point $x\notin\omega$.

In this case, the total actions of the two paths can be improved by switching to paths with vertices
$p,p_1,\ldots,p_j,q_{k+1},\ldots,q_m$ and, respectively, $p,q_1,\ldots,q_{k},p_{j+1},\ldots,p_n.$
%
Therefore, we obtain a contradiction with the optimality of the original paths.
\epf

\begin{lem}\label{lem:uniqueness-of-geodesic}
Let $v\in \R$. With probability one, every Poisson point belongs to at most one semi-infinite minimizer
with asymptotic slope $v$.
\end{lem}
\bpf Let ${\cal U}(v)$ be the event that two distinct
semi-infinite minimizers with asymptotic speed $v$ pass through a common Poissonian point $p$.
Let
\[{\cal S}=\{v\in \R\ :\ \P({\cal U}(v))>0\}.\]
It is sufficient to  show that ${\cal S}=\emptyset$. The invariance under shear transformations implies that either
${\cal S}=\emptyset$ or ${\cal S}=\R$. The latter will be excluded as soon as we prove that ${\cal S}$ is at
most countable.

A triple of
distinct Poisson points $(p,q_1,q_2)$ is called a bifurcation triple for $v$ if there exist two distinct semi-infinite
minimizers $\gamma_1$ and $\gamma_2$ with asymptotic slope $v$ that both start at $p$, then one goes directly (at
constant velocity) to $q_1$ and the other goes directly to $q_2$. We choose $q_1$ such that $\gamma_1$ lies to the left
of $\gamma_2$.

Lemma~\ref{lem:geodesics-do-not-share-more-than-one} implies that  $\gamma_1$ and $\gamma_2$ will not
meet again after $p$ with probability one.

Clearly, ${\cal U}(v)$ implies the existence of such a bifurcation triple for $v$. If $\P({\cal U}(v))>0$, then there
exists $L=L(v)\in {\mathbb N}$ such that the event
\[ {\cal T}_L(v) = \{ \exists \mbox{ bifurcation triple}\ (p,q_1,q_2)\ \mbox{for } v\mbox{ inside the cube }
K((0,0),L)\}\]
has positive probability (using translation invariance).

Now suppose that ${\cal S}$ contains an uncountable number of asymptotic slopes. This implies that there exist
$m,L\in {\mathbb N}$ such that for uncountably many $v$ we would have that
\begin{equation}\label{eq:biftriple}
 \P({\cal T}_L(v))> 1/m.
\end{equation}
Now note that three Poisson  points can only form a bifurcation triple for one $v\in \R$, since otherwise two distinct
semi-infinite minimizers, both starting at some $p$, would have to split up at $p$ and cross again at a later time,
which contradicts Lemma~\ref{lem:geodesics-do-not-share-more-than-one}. Suppose $v_1,v_2,\ldots$ satisfy
\eqref{eq:biftriple}. Denote by $N_L$ the number of Poisson points in the cube $K((0,0),L)$. Then
\[ \sum_{n\geq 1} 1_{{\cal T}_L(v_n)} \leq N_L^3.\]
Taking the expectation on both sides and using \eqref{eq:biftriple} leads to a contradiction.

Therefore, there can be only countably many $v$'s in $\cal S$ which completes the proof.
\epf

With uniqueness in hand we can strengthen Lemma~\ref{lem:exist-geodesic}:

\begin{lem}\label{lem:convergence-of-geodesics}
 With probability one, for every $v\in\R$ and for every sequence $(y_n,t_n)\in \R^2$ with $t_n\to \infty$ and
\[ \lim_{n\to \infty} \frac{y_n}{t_n} = v,\]
and for every Poissonian point $p\in \R^2$,  the minimizing paths
$\gamma_{p,(y_{n},t_{n})}$ converge to a unique semi-infinite minimizer $\gamma_{p,v}$ starting
at $p$ and with asymptotic speed equal to $v$.
\end{lem}
\bpf Let us assume that the convergence does not hold, i.e., there is a sequence $(n')$
such that the restrictions of $\gamma_{p,(y_{n'},t_{n'})}$ and  $\gamma_{p,v}$ on some finite time interval
$I$ do not coincide for all $n'$. Lemma~\ref{lem:exist-geodesic} allows to choose a subsequence $(n'')$ from
$(n')$ such that for sufficiently large $n''$ the restrictions of $\gamma_{p,(y_{n''},t_{n''})}$ on $I$ coincide with
the restrictions of some infinite one-sided geodesics $\gamma'$. The uniqueness established in
Lemma~\ref{lem:uniqueness-of-geodesic} guarantees that $\gamma'$ coincides with $\gamma_{p,v}$, and the resulting
contradiction shows that our assumption was false, completing the proof.
\epf

\subsection{Coalescence of minimizers} Here we prove that any two one-sided minimizers with the same asymptotic slope
coalesce.

\begin{lem}\label{lem:minimizersdontcross}
With probability one it holds that for every $v\in \R$ and for every pair of semi-infinite minimizers, starting 
at different Poisson points, with asymptotic speed $v$, these minimizers either do not touch, or they coalesce at some
Poisson point.
\end{lem}
\bpf Suppose for some $v\in \R$, there do exist two semi-infinite minimizers with asymptotic speed $v$ that touch, but
do not coalesce. If the two minimizers $\gamma_1$ and $\gamma_2$ contain the same Poisson point $p$, then they must stay
together for all times above $p$ according to Lemma~\ref{lem:uniqueness-of-geodesic}. Therefore, the only option
is that $\gamma_1$ and $\gamma_2$ cross, i.e., they consecutively visit Poissonian points
$p_1,p_2,\ldots$, and, respectively, $q_1,q_2,\ldots$, and $[p_1,p_2]\cap [q_1,q_2] = \{x\}$, for some
$x\in \R^2$.

\begin{figure}
\includegraphics[width=3cm]{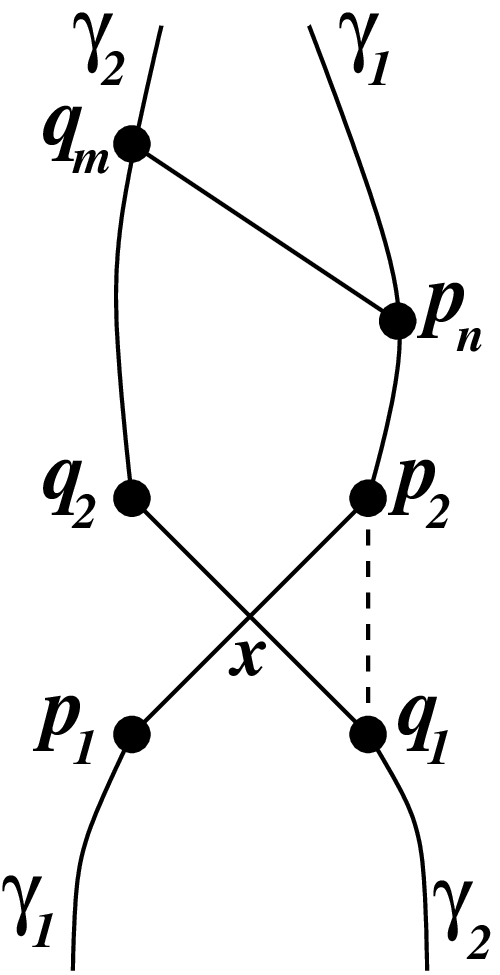} 
\caption{Proof of Lemma~\ref{lem:minimizersdontcross}.}
\label{fig:dont-cross}
\end{figure}

The sequence $\{q_m:m\geq 1\}$ satisfies the conditions of Lemma~\ref{lem:convergence-of-geodesics}, which means
that  the minimizers $\gamma_{p_1,q_{m}}$ converge to $\gamma_1$. However, we claim that with probability 1, none of the
minimizers $\gamma_{p_1,q_m}$ contain any of the $p_n$ $(n\geq 2)$.
In fact, if this claim is violated for some $m,n$, then due to a.s.-uniqueness of a
geodesic between any two Poisson points, we know that $\gamma_{p_1,q_m}$ passes through $p_2$ and $x$. This implies
that action picked up by $\gamma_2$ between $x$ and $q_m$ must be equal to the action picked up by $\gamma_{p_1,q_m}$
between $x$ and $q_m$. However, this contradicts the optimality of $\gamma_2$ as the comparison with the path connecting
$q_1$ directly to $p_2$ and then following $\gamma_{p_1,q_m}$ shows, see Figure~\ref{fig:dont-cross}. The proof is
complete.\epf

For every $v\in\R$, coalescence of one-sided geodesics with asymptotic slope~$v$ generates an equivalence relation on
Poissonian points. We call each equivalence class a coalescence component.

\begin{lem}\label{lem:unbounded_trees} Let $v\in\R$. With probability 1, every coalescence component is unbounded below
in time.
\end{lem}
\bpf Due to shear invariance, it is sufficient to  consider only $v=0$. Suppose that the probability of existence of a
coalescence component bounded from below is
positive. Then there is $(i,j)\in\Z^2$ such that with positive probability the earliest
(i.e., having the minimal time coordinate) point in some coalescence component belongs to $[i,i+1]\times[j,j+1]$. Due to
stationarity, this probability is positive for
any $(i,j)\in\Z^2$. The Tempelman multi-parametric ergodic theorem (see, e.g., \cite[Chapter
6]{Krengel:MR797411})
 implies that there is a constant $c>0$ and a family of random variables $(n_r)_{r\in\NN}$ such that with probability~1,
\begin{equation}
\label{eq:quadratic_growth_of_number_of_trees}
\sum_{i,j\in \Delta(r,n)}I_{ij}> cn^2,\quad n\ge n_r.
\end{equation}
where $\Delta(r,n)=\{(i,j):\ r\le j \le n,\ -j\le i \le j\}$ and $I_{ij}$ is the indicator of event that there is a
tree with earliest point within $[i,i+1]\times[j,j+1]$.

Let $\gamma_+=\gamma_{3,0}$ and $\gamma_-=\gamma_{-3,0}$ be the one-sided optimal paths emitted from
the origin $0$ with asymptotic slopes $3$ and $-3$, respectively. Then there is a random variable $m$ such that
for all $t>m$, one has $-4t<\gamma_-(t)< -2t$ and $2t<\gamma_+(t)<4t$. In particular, the set
\[
\bigcup_{n\ge m}  \bigcup_{(i,j)\in \Delta(m,n)} [i,i+1]\times[j,j+1]
\]
is bounded on the left by $\gamma_-$ and on the right by $\gamma_+$.

Let us denote the coalescence component giving rise to $I_{ij}=1$ by $C_{ij}$. If there are several such components we choose the one
containing the
earliest point. All the components $C_{ij}$, $j>m$ are disjoint and contained in the set
$\{(x,t): t>m,\ -2t\le x\le 2t\}$.

Combining this with \eqref{eq:quadratic_growth_of_number_of_trees}, we can conclude that for every $n\ge n_m$,
the segment
$J_n=[\gamma_-(n),\gamma_+(n)]\times \{n\}$
is crossed by geodesics from at least $cn^2$ disjoint coalescence components. Each of these components has to have an
edge connecting two Poissonian points
and crossing $J_n$.

Therefore, either at least $cn^2$ Poissonian points are contained in $R_n=[-4n,4n]\times[n,n-\sqrt{n}]$, or
there is an edge connecting two vertices of one component and passing through a point on
$H_n=[-4n,4n]\times\{n-\sqrt{n}\}$ and a point
on $L_n=[-4n,4n]\times\{n\}$. Let us denote the former event by $D_n$ and the latter by $E_n$. We have
\begin{align*}
 \P(D_n)\le e^{-cn^2}\E e^{\omega(R_n)}\le e^{-cn^2} e^{8n^{3/2}(e-1)},
\end{align*}
and the Borel--Cantelli Lemma implies that with probability 1 only finitely many events $D_n$ occur.
To prove an analogous statement for the events $E_n$, we need the following lemma:

\begin{lem}\label{lem:no_points_in_parallelogram} Let $x,y\in\R$. Suppose there is an optimal path $\gamma$ connecting a
point $x'\in [x,x+n^{1/5}]\times\{n-\sqrt{n}\}$ straight to a point $y'\in[y,y+n^{1/5}]\times\{n\}$, avoiding all
Poissonian points
between $n-\sqrt{n}$ and $n$. Then the parallelogram $\Pi_n(x,y)$ that is obtained by intersection of
$\R\times [n-\sqrt{n}-1,n-\sqrt{n}+1]$ and the parallelogram with vertices $x,x+n^{1/5},y+n^{1/5},y$, does not
contain
any Poissonian points.
\end{lem}
\bpf Suppose that   $\Pi_n(x,y)$ contains a Poissonian point $p=(z,s)$. Then the original straight path is not optimal,
which follows from the comparison with the path $\tilde \gamma$ consisting of two segments connecting $x'$ to $p$ and
$p$ to $y'$. In fact, an elementary calculation (see also the proof of Lemma~6.5 in \cite{Bakhtin-quasicompact}) shows
that
\[
 S^{n-\sqrt{n},n}(\tilde\gamma)-S^{n-\sqrt{n},n}(\gamma)=\frac{\sqrt{n}}{2(s-n+\sqrt{n})(n-s)}\delta^2,
\]
where $\delta$ is the distance from $p$ to $\gamma$ measured along the spatial axis. Noticing that $\delta\le n^{1/5}$
and $s\in[n-\sqrt{n}-1,n-\sqrt{n}+1]$, we obtain that for large~$n$, the right-hand side is less than $1$, so
increasing the number of points that the path passes through by $1$ overweighs the increase of kinetic
action.
\epf

\smallskip
To finish the proof of Lemma~\ref{lem:unbounded_trees}, we notice that each $H_n$ and $L_n$ can be covered by
$8n^{1-1/5}+1$ intervals of length $n^{1/5}$. The probability of $E_n$ can be then bounded by the sum over all pairs of
these intervals such that there is a straight line minimizer connecting points from these two intervals.
Lemma~\ref{lem:no_points_in_parallelogram} implies that
\begin{align*}
 \P(E_n)&\le (8n^{1-1/5}+1)^2\P\{\omega(\Pi_n(0,0))=0\}\\
        &\le (8n^{1-1/5}+1)^2 e^{-2n^{1/5}}.
\end{align*}
since the area of $\Pi_n(x,y)$ equals $2n^{1/5}$ and does not depend on $x,y$. The
Borel--Cantelli Lemma implies that with probability~$1$, only finitely many events $E_n$ can happen. This completes the
proof of  Lemma~\ref{lem:unbounded_trees}.
\epf

\begin{lem}\label{lem:coalescense of geodesics}
Let  $v\in\R$. With probability one, every two semi-infinite minimizers with asymptotic slope $v$
coalesce.
\end{lem}
\bpf
Due to shear invariance, it is enough to prove the lemma for $v=0$. Also, it is enough to prove it
for minimizers starting at Poisson points. Lemma~\ref{lem:minimizersdontcross} shows that if two of these
semi-infinite minimizers do not coalesce, then they must be disjoint. We call the event that this happens $E$, and
suppose $\P(E)>0$. Now we define $E_k$ $(k\in \mathbb Z)$ as the event that there exist two disjoint minimizers with
asymptotic slope $0$ that both start at Poisson points with times smaller than $k$. Clearly, there exists $k\in \Z$ with
$\P(E_k)>0$, and using translation invariance, we conclude that $\P(E_{-1})>0$. Let us now define, for $x_1,x_2\in
\R\times (-\infty,-1)$ and $\delta\in(0,1)$, $E(x_1,x_2;\delta)$ as the event that there exist two disjoint minimizers
$\gamma_1$ and $\gamma_2$ with asymptotic slope $0$, such that $\gamma_1$ starts within distance $\delta>0$ from $x_1$
and $\gamma_2$ starts within distance $\delta$ from $x_2$. Again it follows that there exist $x_1,x_2\in \R\times
(-\infty,-1)$ and $\delta\in(0,1)$ such that $\P(E(x_1,x_2;\delta))>0$. Let $W_i$ be the crossing of the minimizer
starting near $x_i$ with the $x$-axis. Clearly, there exists $K>0$
such that the event $E(x_1,x_2;\delta,K)$, consisting of all $\omega\in E(x_1,x_2;\delta)$ with $|W_i-x_{i1}|<K$,
has positive probability. Now we choose $h>|x_1-x_2|+2\delta+2K$ and define $z=(h,0)$. Using the
ergodicity of the Poisson point field with respect to translations, we obtain that there are $m,n\in \Z$, $m\neq n$ such
that
\[\P(E(x_1+mz,x_2+mz;\delta,K)\cap E(x_1+nz,x_2+nz;\delta,K))>0.\]
Again using translation invariance, we can take $m=0$ and $n>0$. Let $x_3=x_1+nz$ and $x_4=x_2+nz$.
Because of our choice of $h$, we know that $\max(W_1,W_2)<\min(W_3,W_4)$. Now relabel (if necessary) to
ensure that $W_1<W_2<W_3<W_4$. We denote the minimizer starting near $x_i$ by $\gamma_i$. By construction,
$\gamma_1$ and $\gamma_2$ do not cross, nor do $\gamma_3$ and $\gamma_4$. This implies that the three minimizers
$\gamma_1,\gamma_2$ and $\gamma_4$ do not cross for times larger than $0$, since $\gamma_3$ lies between $\gamma_2$ and
$\gamma_4$ (in principle, $\gamma_2$ and $\gamma_3$ could still coalesce). Our conclusion is that there is a positive
probability of having three non-intersecting minimizers with asymptotic speed $0$, all three starting below time $0$.\\

We can now conclude that for $\delta>0$ small enough and $K,N>0$ big enough, there exist $R>0$, $x\in \R\times
\{-\delta\}$, $y_1,\ldots y_n \in \R\times (-\infty, 0)$ and $z_1,\ldots ,z_m\in\R\times (-\infty, 0)$ with $y_{n,2}\leq
-K-\delta$ and $z_{m,2}\leq -K-\delta$, such that with positive probability:
\begin{itemize}
\item There exist three disjoint semi-infinite minimizers $\gamma_1,\gamma_2$ and $\gamma_3$ with asymptotic slope
$0$ and $\gamma_1(0)<\gamma_2(0)<\gamma_3(0)$, where $\gamma_1$ starts at a Poisson point $p_1$ in $y_n+[0,\delta]^2$,
$\gamma_3$ starts at a Poisson point $p_3$ in $z_m+[0,\delta]^2$ and $\gamma_2$ starts at a Poisson point $p_2$ in
$x+[0,\delta]^2$.
\item If
$p_{1,2}<p_{3,2}$, then $|p_{3,1}-\gamma_1(p_{3,2})|<N$; if  $p_{1,2}\ge p_{3,2}$, then $|p_{1,1}-\gamma_3(p_{1,2})|<N$.
\item There is exactly one Poisson point in $x+[0,\delta]^2$ and the next Poisson point of $\gamma_2$ (above $x$) has a
strictly positive time.
\item There is exactly one Poisson point in $y_i+[0,\delta]^2$ and in $z_j+[0,\delta]^2$, for all $1\leq i\leq n$ and $1\leq j\leq m$.
\item All Poisson points of $\gamma_1$ with negative time lie in $\cup_{i=1}^n( y_i+[0,\delta]^2)$.
\item All Poisson points of $\gamma_3$ with negative time lie in $\cup_{j=1}^m (z_j+[0,\delta]^2)$.
\item $\gamma_1$ and $\gamma_3$ cross the $x$-axis in the interval $[-R,R]$.
\end{itemize}
Let us call this event $A_{\delta,K,N}$. We choose $M>R$ and $L>K+\delta$ such that $[-M,M-\delta]\times [-L,0]$
contains all $y_i$'s, all $z_j$'s and $x$. If $\omega\in A_{\delta,K,N}$, then we can modify the Poisson configuration
$\omega$ by deleting all Poisson points in
\[ [-M,M]\times [-L,0] \setminus \left(\bigcup_{i=1}^n (y_i+[0,\delta]^2)\,\cup\, \bigcup_{j=1}^m
(z_j+[0,\delta]^2)\, \cup\, (x+[0,\delta]^2)\right),\]
and call this new configuration $\tilde{\omega}$. Clearly, $\gamma_1, \gamma_2$ and $\gamma_3$ are still semi-infinite
minimizers under $\tilde{\omega}$. Furthermore, Lemma 8 of \cite{HoNe3} shows that there exists a measurable event
$\tilde{A}_{\delta,K,N}$ with positive probability, such that
\[ \tilde{A}_{\delta,K,N} \subset \{ \tilde{\omega}\ :\ \omega \in A_{\delta,K,N}\}.\]

We claim that on the event $\tilde{A}_{\delta,K,N}$, any semi-infinite minimizer $\gamma$, starting at
 a Poisson point $p$ in the set $(\R\times (-\infty,0] )\setminus ([-M,M]\times [-L,0])$ cannot coalesce with
$\gamma_2$. Therefore, the coalescence component of $\gamma_2$ is bounded below by $-L$ which contradicts
Lemma~\ref{lem:unbounded_trees} and finishes the proof.

Suppose our claim is wrong and $\gamma$ connects to $\gamma_2$. Since $\gamma$ cannot cross $\gamma_1$ or $\gamma_3$, this implies that $\gamma$ passes \emph{between} $\gamma_1$ and $\gamma_3$, and thus at most at distance $N$ from either $p_1$ or $p_3$, when going with constant
velocity from the boundary of $[-M,M]\times [-L,0]$ to $p_2$ or a point on the $x$-axis between $\gamma_1$ and $\gamma_3$.
Therefore, fixing $N$, choosing
$K$ large enough, and increasing $M$ and $L$ if necessary, we can guarantee that if we modify $\gamma$ and let it
pick up one additional Poissonian point $p_1$ or $p_3$, we will decrease its action which  contradicts the optimality
of $\gamma$. The proof is completed. \epf

\subsection{Domains of influence of Poisson points}
Let us fix a mean velocity value $v\in\R$ for the rest of this section. By a one-sided optimal path we will always mean
a one-sided optimal path of asymptotic slope $v$. In this section we describe the random partition of the real line
$R\times\{0\}$ (or, equivalently, $R\times\{t\}$ for any $t\in\R$) into domains of influence of Poisson points. The
domain of influence of $p\in\omega$ consists of all points $x\in\R$ such that the first Poisson point
visited by a unique one-sided minimizer with endpoint $(x,0)$ is $p$. With a little more effort one can give a
description of a space-time tessellation of $\R^2$ of the same kind, but we
omit it for brevity.

\begin{lem}\label{lem:limit-of-minimizers} With probability one the following holds true. For any
convergent sequence of real numbers $x_n\to x_\infty$, and any sequence $\gamma_{n}$ of one-sided
minimizers with endpoints $(x_n,0)$, if there is a point $q$ such that the first Poissonian point visited by each
$\gamma_{n}$ is $q$, then (i) after $q$ all these paths are  uniquely defined and coincide, and (ii) the path
$\gamma_{\infty}$ starting at $(x_\infty,0)$ going straight to $q$ and coinciding with paths $\gamma_{n}$
after $q$, is also a one-sided optimal path.
\end{lem}
\bpf After $q$, the one-sided optimal path is uniquely defined due
to Lemma~\ref{lem:uniqueness-of-geodesic}. Assuming that there is a finite improvement of $\gamma_\infty$, i.e.,
another path starting at $(x_\infty,0)$ connecting to $\gamma_1$ above $q$, one can also construct an improvement of
a path $\gamma_n$ for sufficiently large $n$, i.e., for $(x_n,0)$ sufficiently close to $(x_\infty,0)$, a contradiction.
\epf

\begin{lem}\label{lem:finitely-many-Poissonian-directions} With probability one, for every bounded set $S\in\R^2$
there is a finite set $P$ of Poissonian points such that for every one-sided optimal path
originating at any $p\in S$, the first Poissonian point it visits belongs to $P$.
\end{lem}
\bpf
It is easily verified that for any $v>0$ there there is a random $t_0>0$ such that no Poissonian point $(x,t)$
satisfying $|x|\le vt$ and $t\ge t_0$ can be the first point visited by a one-sided optimal path originating within
$K((0,0),1)$.
Also, for any $R,T>0$, there are only finitely many Poissonian points $(x,t)$ satisfying $|x|\le R$ and $t<T$.
Therefore, if there are countably many one-sided minimizers originating from a point within $K((0,0),1)$, then for any
$R>0$ and any $v>0$ there is a Poissonian point $(x,t)$ satisfying $|x|>\max\{R,vt\}$ visited first by a
one-sided optimal path originating in $K((0,0),1)$.

Since one-sided minimizers cannot cross each other, we see that
on the above event, there is a point $p=(y,s)\in K((0,0),1) $ such that one of the rays $[y,+\infty)\times \{s\}$,  or
$(-\infty,y]\times \{s\}$ cannot be crossed by one-sided minimizers originating at negative times.  So, if this happens
with positive probability, then, due to stationarity and ergodicity, with probability one, there are infinitely many
Poissonian points $(y_n,s_n)_{n\in\Z}$
satisfying $s_n\in[0,1]$ and $\lim_{n\to\pm\infty} y_n=\pm\infty$ such that one of the rays $[y_n,+\infty)\times
\{s_n\}$, or $(-\infty,y_n]\times \{s_n\}$ cannot be crossed by any one-sided optimal path starting at a negative time. Therefore, no one-sided
optimal path can originate at a negative time. This contradiction proves the lemma for $S=K((0,0),1)$. The lemma in
full generality follows due to stationarity, since one can cover any bounded $S$ by finitely many
translates of $K((0,0),1)$.
\epf

\begin{lem}\label{lem:continuity-and-uniqueness-of-geodesics} The following holds with probability one for all $x\in\R$
simultaneously. Suppose there are $n\in\NN$
one-sided optimal paths originating at $(x,0)$. Let us denote by $(x_1,t_1)$,\ldots,$(x_n,t_n)$ the first Poissonian
points visited by these paths, ordering them so that $(x_1-x)/t_1<\ldots<(x_n-x)/t_n$ (note that we still allow $n=1$). Then there is $\eps>0$ such that
(i) for all $y\in(x-\eps,x)$, there is a unique one-sided optimal path originating at
$(y,0)$, and the first Poissonian point visited by it is $(x_1,t_1)$, and (ii) for all $y\in(x,x+\eps)$, there is a unique one-sided optimal path originating at $(y,0)$, and the first Poissonian point visited by it is $(x_n,t_n)$.
\end{lem}
\bpf We will prove only part (ii), since the proof of part (i) is the same.

Suppose that there is a sequence $y_k\downarrow x$ and a sequence $(\gamma_k)$ of one-sided
minimizers originating at $y_k$ with first visited Poissonian point $q_k\ne (x_n,t_n)$.
Lemma~\ref{lem:finitely-many-Poissonian-directions} implies
that the set $\{q_k\}$ is finite, and we can choose a Poissonian point $q=(x',t')$ from this set and a subsequence
$y_{k'}$ such that the first Poissonian point visited by the corresponding paths $\gamma_{k'}$ is $q$. According to
Lemma~\ref{lem:limit-of-minimizers}, there is an optimal path connecting $x=\lim {y_{k'}}$ to $q$. However, optimal
paths cannot cross, and it follows from our construction that $(x'-x)/t'>(x_n-x)/t_n$, and thus $q$ cannot belong to
the set of $n$ Poissonian points in the statement of the Theorem. The resulting contradiction finishes the proof.
\epf

\begin{lem}\label{lem:domains-on-real-line} With probability 1, there is an increasing doubly infinite sequence of
points $(x_k)_{k\in\Z}$ on the
 real line $\R$ satisfying $\lim_{k\to\pm\infty}x_k=\pm\infty$ with the following properties:
\begin{enumerate}
 \item For every $k\in\Z$ there is a Poisson point $p_k$ such that for all $x\in (x_k,x_{k+1})$, there is a unique one-sided
minimizer $\gamma_{(x,0),v}$ originating at $(x,0)$ and the first Poissonian point it visits is $p_k$.
\item For every  $k\in\Z$ there are at least two one-sided optimal paths originating at $x_k$, they pass through $p_k$
and $p_{k-1}$, respectively.
\item For any $k\in\Z$ any $x\in(x_k,x_{k+1})$ and any sequence $(y_n,t_n)\in \R^2$ with $t_n\to \infty$
and
\[ \lim_{n\to \infty} \frac{y_n}{t_n} = v,\]
and for every  $x\in (x_k,x_{k+1})$,  the minimizing paths
$\gamma_{(x,0),(y_{n},t_{n})}$ converge to $\gamma_{(x,0),v}$.
\end{enumerate}
\end{lem}
\bpf The first two parts follow from Lemma~\ref{lem:continuity-and-uniqueness-of-geodesics}. In particular,
that lemma shows that each $x\in\R$ has a neighborhood of points
(excluding $x$) with a unique minimizer. The set of all $x\in\R$ such that $(x,0)$ has two or more one-sided minimizers
must therefore be discrete (it has no accumulation points); this set will be $\{x_k:\ k\in \Z\}$. The proof of the
last part repeats that of Lemma~\ref{lem:convergence-of-geodesics}.
\epf

\section{Busemann functions and stationary solutions of the  Burgers equation}\label{sec:global_solutions}

In this Section we use the one-sided minimizers to construct global solutions of the Burgers equation, thus proving the
existence part of Theorem~\ref{thm:global_solutions}.

Let us summarize some facts on one-sided backward minimizers that follow from Section~\ref{sec:one-sided_minimizers}.
For any velocity $v\in\R$, the following holds with
probability 1.
For every point $p=(x,t)$ there is a non-empty set $\Gamma_{v,p}$ of one-sided action minimizers
$\gamma:(-\infty,t]\to\R$ with asymptotic slope $v$
\[
 \lim_{s\to-\infty}\frac{\gamma(s)}{s}=v,
\]
ending at $p$.
They all coalesce, i.e., they coincide on $(-\infty,t_{v,p}]$ for
some $t_{v,p}<t$. For most points~$p\in\R^2$, $\Gamma_{v,p}$ consists of a unique minimizer $\gamma_{v,p}$, but even if
the uniqueness does not hold, there is the right-most minimizer $\gamma_{v,p}\in\Gamma_{v,p}$ such that
$\gamma_{v,p}(s)\ge \gamma(s)$ for $s\le t$ and any other minimizer $\gamma\in\Gamma_{v,p}$.

For every two points $p_1=(x_1,t_1)$ and $p_2=(x_2,t_2)$, all their one-sided minimizers coalesce, i.e.,
there is
a time $t_v=t_v(p_1,p_2)$ such that $\gamma_{v,p_1}(s)=\gamma_{v,p_2}(s)$ for all $s\le t_v$.

This allows us to define Busemann functions for slope $v$:
\[
 B_v(p_1,p_2)= B_{v,\omega}(p_1,p_2)= A_\omega^{t_v(p_1,p_2),t_2}(\gamma_{v,p_2}) -
A_\omega^{t_v(p_1,p_2),t_1}(\gamma_{v,p_1}),\ p_1,p_2\in\R^2.
\]
Although $t_v$ is not defined uniquely, the definition clearly does not depend on a
concrete choice of $t_v$ or $\gamma_{v,p_1},\gamma_{v,p_2}$. One can also choose $t_v$ to be the maximal of all possible
coalescence times.

Some properties of Busemann functions are summarized in the following lemma:
\begin{lem} \label{lem:busemann_properties} Let $B_v$ be defined as above for $v\in\R$.
 \begin{enumerate}
  \item The distribution of $B_v$ is translation invariant: for any $\Delta \in\R^2$,
\[
 B_v(\cdot+\Delta,\cdot+\Delta)\stackrel{distr}{=} B_v(\cdot,\cdot).
\]
\item
 $B_v$ is antisymmetric:
\[
 B_v(p_1,p_2)=-B_v(p_2,p_1),\quad p_1,p_2\in\R^2,
\]
in particular $B_v(p,p)=0$ for any $p\in\R^2$.
\item $B_v$ is additive:
\[
 B_v(p_1,p_3)=B_v(p_1,p_2)+B_v(p_2,p_3),\quad  p_1,p_2,p_3\in\R^2.
\]
\item  For any $p_1,p_2\in\R^2$, $\E |B_v(p_1,p_2)|<\infty$. \label{item:expectation_finite}
 \end{enumerate}
\end{lem}
\bpf
The first three parts of the Lemma are straightforward. Let us prove part~\ref{item:expectation_finite}.

Using the additivity and translation invariance of Busemann functions, we see that it is sufficient to consider
points $(0,0)$ and $(x,t)$ with $t< 0$. Since the effect of shear transformations on Busemann functions is easily
computable, it is sufficient to assume that $v=0$.

If $s<t<0$, and $x,y\in\R$ we have
\[
A^{s,0}(y,0)\le A^{s,t}(y,x)+\frac{x^2}{2|t|},
\]
so
\[
B_0((0,0),(x,t))\ge  - \frac{x^2}{2|t|},
\]
and all we need is an upper bound on $\E B_0((0,0),(x,t))$.

Since for any $s<t-1$ and any $y\in\R$ we have
\[
A^{s,t}(y,x) \le A^{s,t-1}(y,0)+\frac{x^2}{2},
\]
it is sufficient to assume $x=0$. Also, if $t\in(-1,0)$, then
\[
 A^{s,t}(y,0)\le A^{s,-1}(y,0),
\]
so it is sufficient to consider $t\le -1$. If we prove finiteness of expectation for $t=-1$, then it
will also follow for $t=-2,-3,\ldots$  by additivity, and for all intermediate times by
\[
 A^{s,[t]+1}(y,0)\le A^{s,t}(y,0)\le A^{s,[t]}(y,0).
\]
Therefore, it remains to prove that
\begin{equation}
\label{eq:need_to_estimate_expectation_from_above}
\E (A^{s,-1}(y,0)-A^{s,0}(y,0))<\infty.
\end{equation}
where $(y,s)$ is the space-time point of coalescence of one-sided minimizers with zero asymptotic slope for the points
$(0,0)$ and $(0,-1)$.

\smallskip
Let $H=\{((x,t): t\le -1-|x|)\}$. Notice that $(0,-1)$ is the vertex of the right angle formed by $\partial H$.

Since we are considering minimizers with zero asymptotic slope, we can define
\[
\tau=\inf\{t\le 0:\ (\gamma^*(t),t)\not\in H\},
\]
where  we introduced $\gamma^*=\gamma_{0,(0,0)}$ for brevity. Let us denote
$z=\gamma^*(\tau)=\pm (\tau+1)$.
 If $s\in[\tau,0]$, then
\[
 A^{s,0}(y,0)-A^{s,-1}(y,0)=A^{\tau,0}(z,0)-A^{\tau,-1}(z,0).
\]
If $s<\tau$, then
\begin{eqnarray*}
A^{s,-1}(y,0)-A^{s,0}(y,0) & \leq & A^{s,\tau}(y,z)+A^{\tau,-1}(z,0)-(A^{s,\tau}(y,z)+A^{\tau,0}(z,0))\\
& \le& A^{\tau,-1}(z,0)-A^{\tau,0}(z,0).
\end{eqnarray*}
Combining last two relations with~\eqref{eq:need_to_estimate_expectation_from_above}, we see that we need to establish
\begin{equation}
\label{eq:finiteness_of_expectation1}
\E A^{\tau,0}(z,0)>-\infty
\end{equation}
and
\begin{equation}
\label{eq:finiteness_of_expectation2}
\E A^{\tau,-1}(z,0)<\infty.
\end{equation}
To prove~\eqref{eq:finiteness_of_expectation1} it is sufficient to show that
\[
\sum_{r=1}^\infty\P\left\{-A^{\tau,0}(z,0)>r\right\}<\infty.
\]
Let us choose a small number $\eps>0$ and write
\begin{equation}
\label{eq:z-far-or-deviation-big}
\P\left\{-A^{\tau,0}(z,0)>r\right\}\le \P\{|z|>\eps r\}+\P\left\{|z|\le\eps r;\ -A^{-1-|z|,0}(z,0) >r\right\}.
\end{equation}

The second term can be estimated via Theorem~\ref{thm:concentration_around_alphat}. Estimating the action of
the motion with unit speed from $(\pm \eps r, -\eps r-1)$ to $(z,\tau)$, we obtain that for $|z|\leq \eps r$,
\[  A^{-1-\eps r,0}(\pm \eps r,0) \leq \frac12 (-|z|-1 - (-\eps r -1)) + A^{-1-|z|,0}(z,0),\]
and therefore
\begin{align*}
 &\P\left\{|z|\le\eps r;\ A^{-1-|z|,0}(z,0)<-r\right\}
\\ \le &
\P\left\{ A^{-1-\eps r,0}(\eps r,0)<-r + \frac{\eps r}{2}
\right\}+ \P\left\{ A^{-1-\eps r,0}(-\eps r,0)<-r + \frac{\eps r}{2}
\right\}\\
= & 2 \P\left\{ A^{\eps r+1} + \frac12 \frac{(\eps r)^2}{1+\eps r} < -r + \frac{\eps r}{2}
\right\}
\\
\le & 2 \P\left\{ A^{\eps r+1}- \alpha(0)(\eps r+1) < -r + \frac{\eps r}{2} - \alpha(0)(\eps r+1)
\right\}
\\
\le & 2 c_1(\eps)\exp\left\{-c_2(\eps)\frac{r^{1/2}}{\ln r}\right\},
\end{align*}
for $r>c_0(\eps)$, where $\eps$ is chosen so small that $1-\eps/2 +\eps\alpha(0)>0$. The right-hand side is summable in
$r$.

Let us now estimate $\P\{z>\eps r\}$, one half of the first term in the r.h.s.\ of~\eqref{eq:z-far-or-deviation-big}.

Let $p_s=(s,-2s)$, $s>0$. Let us find $n_0>0$ such that $\Co(p_s,(2s)^{-\delta})$ does not intersect
$\{(x,-3x):\ x>0\}$ for all $s>n_0$. Since the asymptotic slope of $\gamma^*$ is zero, it will eventually cross
$\{(x,-3x):\ x>0\}$. Therefore, if $z>\eps r>n_0$, then there must be $s>\eps r$ and
a point $q=(q_1,q_2)$ outside of the cone
$\Co(p_s,(2s)^{-\delta})$ such that $\gamma^*$ connects $q$ to $(0,0)$ and passes through
$p_s$.
Using Corollary~\ref{cor:delta-straightness-prob}, we see that if $\eps
r>n_0\wedge M$,
\begin{equation}
\label{eq:tail_of_z}
 \P\{z>\eps r\}\le  c_1e^{-c_2(\eps r)^{\kappa}},
\end{equation}
which is summable in $r$.

This finishes the proof of~\eqref{eq:finiteness_of_expectation1}. To prove~\eqref{eq:finiteness_of_expectation2},
it is sufficient to notice that
\[
 \E A^{\tau,-1}(z,0)\le \E \frac{|z|}{2}
\]
and apply \eqref{eq:tail_of_z}. The proof of part~\ref{item:expectation_finite} of
Lemma~\ref{lem:busemann_properties} is completed. \epf

\medskip

Having the Busemann function at hand, one can define
\[
 U_v(x,t)=B((0,0),(x,t)),\quad (x,t)\in\R^2.
\]

The main claim of this Section is that thus defined $U_v$ is skew invariant under of the
HJHLO cocycle, and its space derivative is the global solution
of the Burgers equation.

Let us recall that the HJHLO evolution is given by
\begin{equation}
\label{eq:Burgers_dynamics_on_potentials}
 \Phi^{s,t}W(y)=\inf_{x\in\R} \{W(x)+A^{s,t}(x,y)\}, \quad s\le t,\quad y\in\R,
\end{equation}
where $A^{s,t}(x,y)$ has been defined in~\eqref{eq:optimal_action_between_two_points}.

\begin{lem}\label{lem:global_solution_of_HJ} Function $U_v$ defined above is a global solution of the Hamilton--Jacobi equation. If $s\le t$, then
\begin{equation*}
 \Phi^{s,t} U_v(\cdot,s) (x)=U_v(x,t).
\end{equation*}
\end{lem}
\bpf Let $\gamma_v$ be a minimizer through $(x,t)$ with slope $v$. Then
\begin{align*}
 U_v(x,t)=&U_v(\gamma_{v}(s),s)+(U_v(x,t)-U_v(\gamma_{v}(s),s))\\
                      =&U_v(\gamma_{v}(s),s) + A^{s,t}(\gamma_v(s),x).
\end{align*}
We need to show that the right-hand side is the infimum of $U_v(y,s)+A^{s,t}(y,x)$ over
all $y\in\R$. Suppose that for some $y\in\R$,
\begin{equation}
 U_v(y,s) + A^{s,t}(y,x) < U_v(\gamma_{v}(s),s) + A^{s,t}(\gamma_v(s),x).
\label{eq:suppose_not_infimum}
\end{equation}

Let us take any minimizer $\bar \gamma_v$ originating at $(y,s)$ and denote by $\tau<s$ the time
of coalescence of $\bar \gamma_v$ and $\gamma_v$. We claim that
\begin{equation}
 A^{\tau,s}(\bar \gamma_v)+A^{s,t}(y,x)< A^{\tau,t}(\gamma_v(\tau),x)= A^{\tau,t}(\gamma_v),
\label{eq:gamma_bar_better}
\end{equation}
which contradicts the minimizing property of $\gamma$. In fact, \eqref{eq:gamma_bar_better} is a consequence of
\begin{align*}
 A^{\tau,s}(\bar \gamma_v)-A^{\tau,s}(\gamma_v)&=U_v(y,s) - U_v(\gamma_{v}(s),s)\\
&< A^{s,t}(\gamma_v(s),x) -A^{s,t}(y,x).
\end{align*}
where the second inequality follows from \eqref{eq:suppose_not_infimum}.
\epf

Another way to approach the Burgers equation is to consider, for $p=(x,t)$,
\[
 u_v(x,t)=\dot\gamma_{v,p}(t).
\]
Then $U_v(x,t)-U_v(0,t)=\int_0^x u_v(y,t)dy$.
We recall that $\Psi^{s,t}w$ denotes the solution at time $t$ of the Burgers equation with initial condition $w$
imposed at time
$s$.
\begin{lem} The function $u_v$ defined above is a global solution of the Burgers equation. If $s\le t$, then
\[
 \Psi^{s,t} u_v (\cdot,s) = u_v(\cdot,t),\quad s\le t.
\]
\end{lem}
\bpf This statement is a direct consequence of Lemmas~\ref{lem:evolution_on_Burgers_potentials},
\ref{lem:global_solution_of_HJ}, and
the definition of the Burgers cocycle $\Psi$.\epf

The function $u_v(\cdot,t)$ is clearly piecewise linear with respect to the space coordinate, with downward jumps,  each
linear regime corresponding to the configuration point visited last by one-sided minimizers, see
Lemma~\ref{lem:domains-on-real-line}.

To prove that $U_v(\cdot,t)\in \HH(v,v)$ for all $t$, we will compute the expectation of its spatial increments
(we already know that it is well defined due to part~\ref{item:expectation_finite} of
Lemma~\ref{lem:busemann_properties}), and prove that $u_v(\cdot,t)$ is mixing with respect to the spatial variable.

\begin{lem}\label{lem:mean_increment} For any  $(x,t)\in\R^2$,
\[
 \E (U_v(x+1,t)-U_v(x,t))=\E B_v((x,t),(x+1,t))=v.
\]
\end{lem}
\bpf
First, we consider the case $v=0$. Due to the distributional invariance of Poisson process under reflections,
\[
\E B_0((x,t),(x+1,t)) = \E B_0((x+1,t),(x,t)).
\]
Combining this with the anti-symmetry of $B_0$, we obtain $\E B_0((x+1,t),(x,t))=0$, as required.

In the general case, we can apply the shear transformation $L:(y,s)\mapsto (y+(t-s)v,s)=(y+vt-vs,s)$. Due to Lemma
\ref{lem:shear},
the one-sided minimizers
of
slope $v$ will be mapped onto one-sided minimizers of slope $0$ for the new Poissonian configuration $L(\omega)$.
We already know that
\[
 \E B_{0,L(\omega)}((x+1,t),(x,t))=0,
\]
and a direct computation based on Lemma~\ref{lem:shear} gives
\[
  B_{0,L( \omega)}((x,t),(x+1,t))=B_{v,\omega}((x,t),(x+1,t))+v,
\]
and our statement follows since $L$ preserves the distribution of Poisson process.
\epf

\begin{lem}\label{lem:mixing} Let $v\in\R$. For any $t$, the process $u_v(\cdot,t)$ is mixing.
\end{lem}
\bpf  Due to translation invariance it is sufficient to consider $t=0$. Using the shear invariance we can restrict
ourselves to the case $v=0$. Notice that the values of $u_v(\cdot,0)$ on an interval
are
determined by
the increments of $U_v(\cdot,0)$ on that interval. Therefore, in this proof we can work with these increments.

Let us fix $a>0$ and $\eps>0$. For any $h>2a$, we consider the
processes
\[L_0(x) = B_0((0,0),(x,0))=U_0(x,0)-U_0(0,0), \quad x\in [-a,a],\]
and
\[
L_h(x) = B_0((h,0),(h+x,0))=U_0(x+h,0)-U_0(h,0),\quad x\in [-a,a].
\]
 Next we define the processes
\[
L^t_0(x) = A_*^{-t,0}(0,0) - A_*^{-t,0}(0,x),\quad  x\in [-a,a],
\]
and
\[
L^t_h(x) = A_*^{-t,0}(h,h) - A_*^{-t,0}(h,x+h),\quad  x\in [-a,a],
\]
where $A_*$ is defined as the usual optimal action $A$, with the restriction that the path cannot cross the
vertical line $\chi$ through $(h/2,0)$. This means that $L^t_0$ and $L^t_h$ are by definition independent, since the
first
depends on the Poisson process left of the vertical line, and $L^t_h$ depends on the Poisson process right of the
vertical line.

Note that if the one-sided minimizers for $(-a,0)$ and $(a,0)$ coalesce above  $-t$, and if
the minimizing paths connecting $(0,-t)$ to points  $(-a,0)$ and $(a,0)$ coincide with the respective one-sided
minimizers above their coalescing point and do not cross $\chi$, then $L_0(x)=L^t_0(x)$ for
all
$x\in [-a,a]$. By choosing $t$ large enough and then $h$ large enough, we can ensure that
 \[ \P\{\forall\ x\in[-a,a]\ :\ L_0(x)=L^t_0(x)\} \ge 1-\eps.\]
Since our set-up is symmetric around $\chi$, we also have
 \[ \P\{\forall\ x\in[-a,a]\ :\ L_h(x)=L^t_h(x)\} \ge 1-\eps.\]
Define
\[ C=C_{t,h} = \{ \forall\ x\in[-a,a]\ :\ L_0(x)=L^t_0(x)\ \mbox{and}\ L_h(x)=L^t_h(x)\}.\]
For events $E,F$ for the process $L_0$, we define $\tau_h(F)$ as the ``translated'' event for $L_h$,
 and $E^t$ and $\tau_h(F)^t$ as the corresponding events for $L^t_0$ and $L^t_h$. Then,
\begin{eqnarray*}
|\P(E\cap \tau_h(F)) - \P(E)\P(F)| &\leq & |\P(E\cap \tau_h(F)\cap C) - \P(E)\P(F)|  + 2\eps\\
& = & |\P(E^t\cap \tau_h(F)^t\cap C) - \P(E)\P(F)|  + 2\eps\\
& \leq &  |\P(E^t\cap \tau_h(F)^t) - \P(E)\P(F)|  + 4\eps\\
& = & |\P(E^t)\P(\tau_h(F)^t) - \P(E)\P(F)|  + 4\eps\\
& \leq & |\P(E^t\cap C)\P(\tau_h(F)^t\cap C) - \P(E)\P(F)|  + 8\eps\\
& = & |\P(E\cap C)\P(\tau_h(F)\cap C) - \P(E)\P(F)|  + 8\eps\\
& \leq & 12\eps,
\end{eqnarray*}
and mixing follows due to the arbitrary choice of $\epsilon$. \epf

Combining Lemmas~\ref{lem:mean_increment} and~\ref{lem:mixing}, we conclude that the Birkhoff space averages of $u_v$
have a well-defined,
deterministic limit $v$, so $U_v\in \HH(v,v)$.

\section{Stationary solutions: uniqueness and basins of attraction}\label{sec:attractor}
In this section we prove Theorem~\ref{thm:pullback_attraction} and the uniqueness part in
Theorem~\ref{thm:global_solutions}

The key step in the proof of Theorem~\ref{thm:pullback_attraction} is the following observation.

\begin{lem}\label{lem:asymptotic_slope_in_pullback_attraction}
Let $t\in\R$ and suppose that an initial condition $W$ satisfies one of the
conditions~\eqref{eq:flux_from_the_left_wins},\eqref{eq:flux_from_the_right_wins},%
\eqref{eq:no_flux_from_infinity}.  With probability one, the following holds true for every $y\in\R$. Let  $y^*(s)$ be a
solution of the optimization problem
~\eqref{eq:Burgers_dynamics_on_potentials}. Then
\[
 \lim_{s\to-\infty}\frac{y^*(s)}{s}=v.
\]
\end{lem}

\bpf[Proof of Theorem~\ref{thm:pullback_attraction}]  Let us take any rectangle $Q=[-R,R]\times [t_0,t_1]$ and
set $t=t_1$. We use Lemma~\ref{lem:domains-on-real-line} to find points $a,b\in\R$ satisfying $a<-R<R<b$, not
coinciding with any of the points $x_k,k\in\Z$ and such that one-sided backward minimizers $\gamma_{(a,t),v}$ and $\gamma_{(b,t),v}$
do not cross $Q$.

Applying
Lemma~\ref{lem:asymptotic_slope_in_pullback_attraction} to $x=a,b$, we see that the corresponding points $a^*(s)$ and
$b^*(s)$ satisfy $a^*(s)/s\to v$ and  $b^*(s)/s\to v$ as $s\to-\infty$.

Let $p=(x_0,\tau_0)$ be the point of coalescence of the one-sided minimizers $\gamma_{(a,t),v}$ and $\gamma_{(b,t),v}$.
We automatically have $\tau_0<t_0$.
Lemma~\ref{lem:domains-on-real-line} then implies that there is $\tau_1<\min\{\tau_0,0\}$ such that for $s<\tau_1$,
the  restrictions of the finite minimizers
connecting $(a^*(s),s)$ to $(a,t)$ and $(b^*(s),s)$ to $(b,t)$ on $[\tau_0,t]$ coincide with the restrictions of
$\gamma_{(a,t),v}$ and $\gamma_{(b,t),v}$ (this also implies that we can choose $a^*(s)=b^*(s)$).

Since $Q$ is trapped between $\gamma_{(a,t),v}$ and $\gamma_{(b,t),v}$, and minimizing paths cannot cross each other,
we conclude that for any
$s<\tau_1$, and any $(x,t)\in Q$, the minimizers connecting $(x^*(s),s)$ to $(x,t)$ (where $x^*$ is a solution of the
optimization problem~\eqref{eq:Burgers_dynamics_on_potentials}) have to pass through $p$. In particular, the slopes of
these minimizers determining the evolution of the Burgers velocity field in $[-R,R]$ throughout $[t_0,t_1]$ do not
change (and coincide with the slopes of one-sided backward minimizers) as long as $s<\tau_1$, which completes the
proof. \epf

\medskip

\bpf[Proof of Lemma~\ref{lem:asymptotic_slope_in_pullback_attraction}] We will only prove the sufficiency of
condition~\eqref{eq:no_flux_from_infinity}. The proof of sufficiency of conditions
\eqref{eq:flux_from_the_left_wins} and \eqref{eq:flux_from_the_right_wins} follows the same lines and we omit it.

Let us also restrict ourselves to $t=0$ for simplicity. The proof does not change for other values of $t$.

 Since $y^*$ is increasing in $y$, it is sufficient to show that the conclusion of the lemma holds with
probability 1 for fixed $y$. The stationarity of Poisson point field implies that we can assume $y=0$.

We must show that for any $\eps>0$ it is extremely unlikely for a path $\gamma$ with $\gamma(0)=0$
and $|\gamma(-r)|> \eps r$ to be optimal if $r$ is large. For definiteness, let us work with paths satisfying
$\gamma(-r)>\eps r$.

For any $\delta_0\in(0,-\alpha(0)/3)$ and for sufficiently large $r$,
\[
 W(0) + A^{-r,0}(0,0)< (\alpha(0)+\delta_0)r.
\]
Let us introduce
\[
Q_{ij}=[\eps j +i,\eps j +i+1]\times [-j-1,-j],\quad i,j\ge 0.
\]
If $(x,-r)\in Q_{ij}$ and $W(x) + A^{-r,0}(x,0) < (\alpha(0)+\delta_0)r$, which would be necessary if $x=0^*(-r)$, then
\[
 \inf_{z\in[\eps j +i,\eps j +i+1]} W(z)+ A^{-j-2,0}(\eps j +i+1,0)< (\alpha(0)+\delta_0)j+\frac{1}{2}.
\]
The condition~\eqref{eq:no_flux_from_infinity} at $+\infty$ implies that there is $j_0$ such that for $j>j_0$ and all
$i$,
\[
 \inf_{z\in[\eps j +i,\eps j +i+1]} W(z)> -(j+i) \delta_0.
\]
so there is $j_1$ such that for $j>j_1$,
\[
  A^{-j-2,0}(\eps j +i+1,0) < (\alpha(0) + 2\delta_0)j +\delta_0i + \frac{1}{2} <
j\left(\alpha(0)+3\delta_0+\delta_0\frac{i}{j}\right).
\]
Let us denote by $B_{ji}$ the event described by this inequality.
Due to the Borel--Cantelli lemma, to show that with probability 1, events $B_{ji}$ can happen only for finitely many values
of $j$, it suffices to
show that for some $\beta>0$ and $c>0$,
\begin{equation}
\label{eq:at_most_linear_growth_BC}
\sum_{j\ge c}\sum_{i\le \beta j}\P(B_{ji})<\infty,
\end{equation}
and
\begin{equation}
\label{eq:super_linear_growth_BC}
\sum_{j\ge c}\sum_{i>\beta j}\P(B_{ji})<\infty.
\end{equation}
 Denoting $\alpha_{ji}=\alpha\left(\frac{\eps j+i+1}{j+2}\right)$, using shear and translation invariance,
we
obtain
\[
\P(B_{ji})
 =
\P\left\{A^{j+2}-\alpha(0)(j+2) <
j\left(\alpha(0)+3\delta_0+\delta_0 \frac{i}{j}- \frac{j+2}{j}\alpha_{ji}\right)\right\}.
\]
If $\delta_0$ was chosen sufficiently small, then, using Lemma~\ref{lem:shape-function},  we can find $j_2$ such
that for all $j>j_2$ and all $i$,
\[
 j\left(\alpha(0)+3\delta_0+\delta_0 \frac{i}{j}- \frac{j+2}{j}\alpha_{ji}\right) <
-(j+2)\left(\frac{\eps^2}{3}+\frac{i^2}{2j^2}\right),
\]
so that
\begin{equation}
\label{eq:B_ji}
 \P(B_{ji})\le \P\left\{A^{j+2}-\alpha(0)(j+2)<-(j+2)\left(\frac{\eps^2}{3}+\frac{i^2}{3j^2}\right)\right\}.
\end{equation}

Now~\eqref{eq:at_most_linear_growth_BC} follows (with any $c\ge j_2$ and with arbitrary choice of $\beta$) from
Theorem~\ref{thm:concentration_around_alphat}.

To prove~\eqref{eq:super_linear_growth_BC} we need an auxiliary lemma. In its statement and proof we use the notation
introduced in Section~\ref{sec:subadd}.
\begin{lem}
\label{lem:tails_of_action}
There are constants $c_1,c_2,c_3,X_0,T_0>0$ such that for $t>T_0$, $x>X_0$,
 \[
  \P\{A^t\le  - xt\}\le c_2 e^{-c_3 xt}.
 \]
\end{lem}

Notice that this lemma directly implies~\eqref{eq:super_linear_growth_BC} and thus
completes the proof of Lemma~\ref{lem:asymptotic_slope_in_pullback_attraction}
 if we choose $c>T_0$
and
$\beta$ satisfying
\[
\frac{\eps^2}{3}+\frac{\beta^2}{3}-\alpha(0)>X_0.
\]
It remains to prove Lemma~\ref{lem:tails_of_action}.

\bpf[Proof of Lemma~\ref{lem:tails_of_action}] It is sufficient to prove the lemma for $t\in\NN$, although the values
of constants may need adjustment for general
$t$.
 Let us take $c_4>0$ and write
\begin{equation}
 \P\{A^t\le -xt\}\le \P\{\#\calA \le c_4 xt,\ A^t\le -xt\}+\P\{\#\calA > c_4 xt\}.
\label{eq:estimating_action_linear_tails}
\end{equation}
To estimate the first term on the r.h.s., we can
use~\eqref{eq:latanibound} and
derive that if $c_4$ is chosen small enough to ensure $y_0c_4<2$, where $y_0$ was introduced before
condition~\eqref{eq:cond_on_y}, then
for some constants $c_5,c_6>0$, any $x>1$, and sufficiently large $t$,
\begin{align*}
 \P\{\#\calA \le c_4 xt,\ A^t\le -xt\}\le \P\{N_{\lceil c_4 xt\rceil}\ge xt\}
                                       \le c_5e^{-c_6 xt}.
\end{align*}
The second term on the r.h.s.\ of~\eqref{eq:estimating_action_linear_tails} can be estimated using
Lemma~\ref{lem:latani}.
If $c_4x\geq R$ and $t$ is sufficiently large, then
\begin{align*}
 \P\{\#\calA > c_4 xt\}&\le \sum_{n\ge c_4xt} \P(E_{n,t}) \\
                      &\le \sum_{n\ge c_4xt}  C_1\exp(-C_2n^2/t)\\
                      &\le C'_1\exp(-C'_2x^2t).
\end{align*}
for some constants $C'_1,C'_2>0$,
which completes the proof.\epf

\medskip

\bpf[Proof of uniqueness in Theorem~\ref{thm:global_solutions}] We will prove that any skew-invariant function $u$
with average velocity $v$ coincides with the global solution $u_v$ at time $0$.

Let us take an arbitrary interval
$I=(a,b)$. Lemma~\ref{lem:asymptotic_slope_in_pullback_attraction} implies that for any $W$ satisfying
$\HH(v,v)$, there
is a time $T_0(a,b,W)\ge 0$ such that if $s<-T_0$, then there is a point $a^*\in\R$ that solves the optimization
problem~\eqref{eq:Burgers_dynamics_on_potentials} for $t=0$ and for all points $y\in I$ at once, and the respective
finite minimizers on $[s,0]$ have the same velocity at time 0 as the infinite one-sided minimizers of asymptotic slope
$v$.

Suppose now that $U_\omega(x,t)=U_{\theta^t\omega}(x)$ is a global solution in $\hat \HH(v,v)$. Then
$T_0(a,b,U_{\theta^t\omega}(\cdot))>0$ is a stationary process. In particular, this means that with probability 1, there
is $R>0$ and a sequence of times $s_n\downarrow -\infty$ such that
$T_0(a,b,U_{\theta^{s_n}\omega}(\cdot))<R$ for all $n\in\NN$. Therefore, there is $n$ such that
$s_n <- T_0(a,b,U_{\theta^{s_n}\omega}(\cdot))$. This and the fact that $U$ at time $0$ is the solution of
problem ~\eqref{eq:Burgers_dynamics_on_potentials} for $t=0$, $s=s_n$ and initial condition
$W=U_{\theta^{s_n}\omega}(\cdot)$, we conclude that $U$ and the global solution $U_v$ coincide
on $I$ at time $0$, and the proof is complete.
\epf

\section{Basics on Burgers with Poissonian forcing}\label{sec:proofs_of_basic_burgers_facts}

\bpf[Proof of Lemma~\ref{lem:full_measure_set_where_dynamics_is_defined}]
Let us take $N\in\NN$, $x\in[-N,N]$, $W\in\HH$, and $M\in\NN$  satisfying
$ W(y)\ge -M(|y|+1)$ for all $y\in \R$ and  $|W(0)|<M$. Let us take any $m,n\in\ZZ$ satisfying $m<n$
and
any $t,s\in\R$ satisfying $m<s<t<n$.

Suppose now that a path $\gamma:[s,t]\to\R$ and a constant $L\in\NN$ satisfy $\gamma(t)=x$ and
$\sup_{r\in[s,t]}|\gamma(r)|\in[L,L+1)$. Let us compare $\gamma$
with the
straight path $\gamma_0(r)=x(r-s)/(t-s)$, $r\in[s,t]$. Assuming that
\[
A_\omega^{s,t}(W,\gamma)\le A_\omega^{s,t}(W,\gamma_0),
\]
we obtain
\begin{align*}
M+\frac{x^2}{2(t-s)}&\ge W(0)+\frac{x^2}{2(t-s)}\\
&\ge A_\omega^{s,t}(W,\gamma_0)\\
&\ge A_\omega^{s,t}(W,\gamma)\\
&\ge -M(L+1+1)+\frac{(x-L)^2}{2(t-s)}-\omega([-(L+1),(L+1)]\times[s,t]).
\end{align*}
Therefore, for sufficiently large $L$,
\[
 \omega([-(L+1),(L+1)]\times[m,n])\ge -M(L+3)+\frac{L^2-2Lx}{2(t-s)}\ge -M(L+3)+\frac{L^2-2LN}{2(n-m)}.
\]
Since the r.h.s.\ is quadratic in $L$, and the l.h.s.\ grows linearly in $L$ with probability~1 due to the strong law
of large numbers, we conclude that
this inequality can be true only for finitely many values of $L\in\NN$. Therefore,
under the imposed restrictions on $x,W,s,t$, there is $L_0=L_0(\omega,m,n)<\infty$ such that  the
variational problem~\eqref{eq:optimization_problem} will not change if supplied with an additional restriction
$\sup_{r\in[s,t]}|\gamma(r)|\le L_0$.

Since with probability 1 there are finitely many Poissonian points in $[-L_0,L_0]\times[m,n]$, it is useful to split
the
variational problem into two parts: optimization over paths that do not pass through any Poisson points and over
paths passing through some configuration points.

The existence of an optimal path among those not passing through any Poisson points is
guaranteed by the theory of unforced Burgers equation $u_t+uu_x=0$.

Also, there are only finitely many broken line paths between Poissonian points, and, for the same reason as above,
for each Poisson point
$p$ in $[-L_0,L_0]\times[s,t)$ there is an optimal path among paths not passing through any other Poisson points
and terminating at $p$.  Gluing these paths together, we see that the
extremum in the optimization problem for paths containing Poissonian points is also attained.

Combining the two cases
leads to our claim holding a.s.\ for all $x,W,s,t$, satisfying constraints specified by integers $M,N,m,n$. The
countable intersection of full measure sets over all $M,N,m,n$ still produces a full measure set.
Its time invariance follows since
if $L_0(\omega,m,n)$ is finite for all $m,n$, then $L_0(\theta^t\omega,m,n)$ is finite for all $m,n$ and all $t\in
\R$. The proof is complete.
\epf

\bpf[Proof of Lemma~\ref{lem:evolution_on_Burgers_potentials}] Part~\ref{it:open-domains-of-influence} follows from the
observation that
the action depends on the path continuously if the family of configuration points visited by the path is fixed.

 In the situation where $\gamma$ does not pass through any
configuration points, Part~\ref{it:two-ways-to-compute-velocity} follows from the theory for unforced Burgers equation,
see e.g., \cite[Section 1.3]{Lions-book:MR667669}. If the minimizer $\gamma$ for $x$ has a
straight line segment connecting a configuration point $p=(y,r)$ to $(x,t)$, then
\begin{equation}
 \frac{d}{dx}(\Phi_\omega^{s,t}W)(x)=\frac{d}{dx}\frac{(x-y)^2}{2(t-r)}=\frac{x-y}{t-r}=\dot\gamma(t),
\label{eq:two-ways-to-compute-velocity}
\end{equation}
and Part~\ref{it:two-ways-to-compute-velocity} is proven for points with unique minimizers. The case of the
boundary points where the minimizers are not unique is considered similarly.

Part~\ref{it:piecewise_linearity} is a direct consequence of Part~\ref{it:two-ways-to-compute-velocity}.
Part~\ref{it:continuity_of_HJ_solution} follows from local boundedness of the slope of minimizers. For minimizers not
passing through any configuration points, this is a consequence of the classical theory of unforced Burgers equation,
and for minimizers passing through some configuration points it follows from~\eqref{eq:two-ways-to-compute-velocity}.
\epf

\bpf[Proof of Lemma~\ref{lem:Hamilton-Jacobi}] Let us take a point $(x,t)\in (\R\times(s,\infty))\setminus\omega$  and a
small neighborhood
$O\subset(\R\times ((s+t)/2,\infty))\setminus\omega$ of $(x,t)$. There is a time $s_0>s$ such that
for all $(x',t')\in O$,
no minimizers realizing $\Phi^{s,t'}W(x')$  pass through any configuration point after
$s_0$.
Lemma~\ref{lem:cocycle} implies now that the restriction of $U$ to $O$ coincides with the entropy solution
of~\eqref{eq:Hamilton-Jacobi} with initial data $\Phi_\omega^{s,s_0}W$ imposed at time $s_0$.
\epf

\bpf[Proof of Lemma~\ref{lem:invariant_spaces}] The Lipschitz property required in the definition of $\HH$ follows from
Lemma~\ref{lem:evolution_on_Burgers_potentials}, so we only have to control the behavior of $\Phi_\omega^{s,t}W(x)$ as
$x\to\pm\infty$. It is easy to see that for our almost sure statement, it is sufficient to construct exceptional sets
for integer
times $s,t$ and the limits in definition of $\HH$ and $\HH(v_-,v_+)$ only along
integer values of~$x$.

 We will prove that if there is a constant $M$ such that if
$W(x)>-M(|x|+1)$ for all $x\in\R$, then for any $\eps>0$ and for sufficiently large integer values of $|n|$,
$\Phi_\omega^{s,t}W(n)> -(M+\eps)|n|$. Suppose, for some $n$ this inequality is violated and for some $L\ge 0$, an
optimal
path $\gamma$ realizing $\Phi_\omega^{s,t}W(n)$ satisfies $\sup_{r\in[s,t]}|\gamma(r)-n|\in[L,L+1]$.
Since
\[
 -(M+\eps))|n|\ge\Phi_\omega^{s,t}W(n) \ge -M(|n-L|+1)+\frac{L^2}{2(t-s)}-\omega([n-L,n+L]\times[s,t]),
\]
we obtain
\[
 \omega([n-L,n+L]\times[s,t])\ge \eps |n|-M(L+1)+\frac{L^2}{2(t-s)}.
\]
The l.h.s.\ has the Poisson distribution with mean $2L(t-s)$, and the Borel--Cantelli Lemma implies that
this can hold only for finitely many values of $n$ and $L$, which proves
the first part of the lemma.

The proof of the second part follows the same lines, and we omit it.
\epf

\bibliographystyle{alpha}
\bibliography{Burgers}

\end{document}